\documentclass{article}

\usepackage{colt11e}
\usepackage[round,comma]{natbib}
\usepackage{amsmath,amsfonts,amssymb}

\newcommand{\ifarXiv}[1]{#1}  
\newcommand{\ifCOLT}[1]{}   

\newcommand{\beq}{\begin{equation}}
\newcommand{\eeq}{\end{equation}}

\newcommand{\beqa}{\begin{eqnarray}}
\newcommand{\eeqa}{\end{eqnarray}}

\newcommand{\beqan}{\begin{eqnarray*}}
\newcommand{\eeqan}{\end{eqnarray*}}

\newcommand{\eqdef}{\stackrel{\rm def}{=}}

\renewcommand{\epsilon}{\varepsilon}
\renewcommand{\leq}{\leqslant}
\renewcommand{\geq}{\geqslant}
\renewcommand{\hat}{\widehat}
\newcommand{\wh}{\widehat}
\renewcommand{\P}{\mathbb{P}}
\newcommand{\E}{\mathbb{E}}
\newcommand{\cA}{\mathcal{A}}
\newcommand{\cB}{\mathcal{B}}
\newcommand{\cC}{\mathcal{C}}

\newcommand{\cE}{\mathcal{E}}
\newcommand{\cS}{\mathcal{S}}
\newcommand{\cK}{\mathcal{K}}
\newcommand{\cX}{\mathcal{X}}
\newcommand{\cP}{\mathcal{P}}
\newcommand{\cL}{\mathcal{L}}
\newcommand{\cF}{\mathcal{F}}
\newcommand{\ind}{\mathbb{I}}
\newcommand{\ol}{\overline}
\newcommand{\mt}{\widetilde{\mu}}
\newcommand{\wt}{\widetilde}
\renewcommand{\d}{\mbox{d}}
\newcommand{\nup}{\kappa}
\newcommand{\eps}{\varepsilon}
\newcommand{\norm}[1][\,\cdot\,]{\ensuremath{\left\Arrowvert #1 \right\Arrowvert}}

\newcommand{\kmin}{\cK_{\inf}}
\newcommand{\andAt}{\ \,\, \mbox{\small and} \ \,\, A_{t+1} = a}

\newlength{\minipagewidth}
\newtheorem{remark}{Remark}

\title{A Finite-Time Analysis of Multi-armed Bandits Problems with Kullback-Leibler Divergences}

\author{Odalric-Ambrym Maillard \\
INRIA Lille Nord-Europe \\
France \\
\texttt{\small odalric.maillard@inria.fr}
\And
R{\'e}mi Munos \\
INRIA Lille Nord-Europe \\
France \\
\texttt{\small remi.munos@inria.fr}
\And
Gilles Stoltz \\
Ecole normale sup{\'e}rieure\thanks{CNRS -- Ecole normale sup{\'e}rieure, Paris -- INRIA, within the project-team CLASSIC}, Paris\\
\& HEC Paris \\
France \\
\texttt{\small gilles.stoltz@ens.fr}
}

\begin{document}

\maketitle

\begin{abstract}
We consider a Kullback-Leibler-based algorithm for the stochastic multi-armed bandit problem in the case of distributions with finite supports
(not necessarily known beforehand), whose asymptotic regret matches the lower bound of \cite{Burnetas96}. Our contribution is to provide a finite-time analysis of this algorithm;
we get bounds whose main terms are smaller than the ones of previously known algorithms with finite-time analyses (like UCB-type algorithms).
\end{abstract}

\section{Introduction}
The \textit{stochastic} multi-armed bandit problem, introduced by \citet{ro52}, formalizes the problem of decision-making under uncertainty, and illustrates the fundamental tradeoff that appears between \textit{exploration}, i.e., making decisions in order to improve the knowledge of the environment, and \textit{exploitation}, i.e., maximizing the payoff. \smallskip

{\bf Setting.} In this paper, we consider a multi-armed bandit problem with \textit{finitely} many arms indexed by $\cA$, for which each arm $a \in \cA$ is associated with an unknown and fixed probability distribution $\nu_a$ over $[0,1]$. The game is \textit{sequential} and goes as follows: at each round $t\geq1$, the player first picks an arm $A_t \in \cA$ and then receives a stochastic payoff $Y_t$ drawn at random according to $\nu_{A_t}$. He only gets to see the payoff $Y_t$.

For each arm $a \in \cA$, we denote by $\mu_a$ the expectation of its associated distribution $\nu_a$
and we let $a^\star$ be any optimal arm, i.e., $\ \displaystyle{a^{\star} \in \mathop{\mathrm{argmax}}_{a \in \cA} \, \mu_a\,.}$ \\
We write $\mu^\star$ as a short-hand notation for the largest expectation $\mu_{a^{\star}}$
and denote the gap of the expected payoff $\mu_a$ of an arm $a \in \cA$ to $\mu^\star$ as
$\Delta_a = \mu^\star - \mu_a$.
In addition, the number of times each arm $a \in \cA$ is pulled between the rounds $1$ and $T$
is referred to as $N_T(a)$,
\[
N_T(a) \eqdef \sum_{t=1}^T \ind_{ \{ A_t = a \} }\,.
\]

The quality of a strategy will be evaluated through the standard notion of \textit{expected regret}, which
we recall now. The expected regret (or simply regret) at round $T \geq 1$ is defined as
\begin{equation}
\label{eq:defregr}
R_T \eqdef \E \! \left[ T \mu^{\star} - \sum_{t=1}^T Y_t \right]
= \E \! \left[ T \mu^{\star} - \sum_{t=1}^T \mu_{A_t} \right]
= \sum_{a \in \cA} \Delta_a \,\, \E \bigl[ N_T(a) \bigr]\,,
\end{equation}
where we used the tower rule for the first equality. Note that the
expectation is with respect to the random draws of the $Y_t$ according to the
$\nu_{A_t}$ and also to the possible auxiliary randomizations that the decision-making strategy
is resorting to.

The regret measures the cumulative loss resulting from pulling sub-optimal arms, and thus quantifies
the amount of exploration required by an algorithm in order to find a best arm, since, as (\ref{eq:defregr}) indicates,
the regret scales with the expected number of pulls of sub-optimal arms.
Since the formulation of the problem by \citet{ro52} the regret has been a popular criterion for assessing the quality of a strategy.
\smallskip

{\bf Known lower bounds.} \citet{LaiRo85} showed that
for some (one-dimensional) parametric classes of distributions, any consistent strategy (i.e.,
any strategy not pulling sub-optimal arms more than in a polynomial number of rounds) will
despite all asymptotically
pull in expectation any sub-optimal arm $a$ at least
\[
\E \bigl[ N_T(a) \bigr] \geq \biggl(\frac{1}{\cK(\nu_a,\nu^\star)} + o(1)\biggr)\log(T)
\]
times,  where $\cK(\nu_a,\nu^\star)$ is the Kullback-Leibler (KL) divergence between $\nu_a$ and $\nu^\star$;
it measures how close distributions $\nu_a$ and $\nu^\star$ are from a theoretical information perspective.

Later, \citet{Burnetas96} extended this result to some classes of multi-dimensional parametric distributions
and proved the following generic lower bound: for a given family $\cP$ of possible distributions over the arms,
\[
\E \bigl[ N_T(a) \bigr] \geq
\biggl( \frac{1}{\cK_{\inf}(\nu_a, \mu^\star)} + o(1) \biggr) \log(T)\,,
\qquad \mbox{where} \quad \cK_{\inf}(\nu_a, \mu^\star) \eqdef \inf_{\nu \in\cP: \, E(\nu)> \mu^*} \cK(\nu_a,\nu) \,,
\]
with the notation $E(\nu)$ for the expectation of a distribution $\nu$.
The intuition behind this improvement is to be related to the goal that we want to achieve in bandit problems;
it is not detecting whether a distribution is optimal or not (for this goal, the relevant quantity would be $\cK(\nu_a,\nu^\star)$),
but rather achieving the optimal rate of reward $\mu^\star$ (i.e., one needs to measure
how close $\nu_a$ is to any distribution $\nu\in\cP$ whose expectation is at least $\mu^\star$). \smallskip

{\bf Known upper bounds.} \citet{LaiRo85} provided an algorithm based on the KL divergence, which has been extended by \citet{Burnetas96} to an algorithm based on $\cK_{\inf}$; it is asymptotically optimal since the number of pulls of any sub-optimal arm $a$ satisfies
\[
\E \bigl[ N_T(a) \bigr] \leq \biggl( \frac{1}{\cK_{\inf}(\nu_a, \mu^\star)} + o(1) \biggr) \log(T)\,.
\]
This result holds for finite-dimensional parametric distributions under some assumptions, e.g.,
the distributions having a finite and known support or belonging to a set of Gaussian distributions with known variance.
Recently \citet{Honda10} extended this asymptotic result to the case of distributions $\cP$ with support in $[0,1]$
and such that $\mu^*<1$; the key ingredient in this case is that $\cK_{\inf}(\nu_a,\mu^\star)$
is equal to $\cK_{\min}(\nu_a, \mu^\star) \eqdef \inf_{\nu \in \cP : E(\nu)\geq \mu^*} \cK(\nu_a,\nu)$. \smallskip

{\bf Motivation.} All the results mentioned above provide asymptotic bounds only.
However, any algorithm is only used for a finite number of rounds and it is thus essential to provide a finite-time analysis of its performance. \citet{Auer02} initiated this work by providing an algorithm (UCB1) based on a Chernoff-Hoeffding bound; it pulls any sub-optimal arm,
till any time $T$, at most $(8/\Delta_a^2) \log T + 1 + \pi^2/3$ times, in expectation.
Although this yields a logarithmic regret, the multiplicative constant depends on the gap $\Delta_a^2 = (\mu^\star - \mu_a)^2$ but not on $\cK_{\inf}(\nu_a, \mu^\star)$, which can be seen to be larger than $\Delta_a^2/2$ by Pinsker's inequality;
that is, this non-asymptotic bound does not have the right dependence in the distributions.
(How much is gained of course depends on the specific families of distributions at hand.)
\citet{AudibertTCS09} provided an algorithm (UCB-V) that takes into account the empirical variance of the arms and exhibited
a strategy such that $\E \bigl[ N_T(a) \bigr] \leq 10 (\sigma_a^2/\Delta_a^2 + 2/\Delta_a) \log T$ for any time $T$
(where $\sigma_a^2$ is the variance of arm $a$); it improves over UCB1 in case of arms with small variance. Other variants include the MOSS algorithm by \cite{Audibert2010} and Improved UCB by \cite{Auer2010}.

However, all these algorithms only rely on one moment (for UCB1) or two moments (for UCB-V) of the empirical distributions of the obtained rewards; they do not fully exploit the empirical distributions. As a consequence, the resulting bounds are expressed in terms of the means $\mu_a$ and variances $\sigma_a^2$ of the sub-optimal arms and not in terms of the quantity $\cK_{\inf}(\nu_a, \mu^\star)$ appearing in the lower bounds.
The numerical experiments reported in \cite{Filipi2010} confirm that these algorithms are less efficient than those based on $\cK_{\inf}$. \smallskip

{\bf Our contribution.} In this paper we analyze a $\kmin$-based algorithm inspired by
the ones studied in \cite{LaiRo85,Burnetas96,Filipi2010}; it indeed
takes into account the full empirical distribution of the observed rewards.
The analysis is performed (with explicit bounds) in the case of Bernoulli distributions over the arms.
Less explicit but finite-time bounds are obtained in the case of finitely supported distributions (whose supports
do not need to be known in advance).
Finally, we pave the way for handling the case of general
finite-dimensional parametric distributions.
These results improve on the ones by \citet{Burnetas96,Honda10} since finite-time bounds (implying their asymptotic results) are obtained;
and on \citet{Auer02,AudibertTCS09} as the dependency of the main term scales with $\kmin(\nu_a, \mu^\star)$.
The proposed $\kmin$-based algorithm is also more natural and more appealing than the one
presented in \citet{Honda10}. \medskip

{\bf Recent related works.}
Since our initial submission of the present paper, we got aware of two papers that tackle problems similar to ours.
First, a revised version of \citet[personal communication]{HondaArxiv} obtains finite-time regret bounds
(with prohibitively large constants) for a \emph{randomized} (less natural) strategy in the case of distributions with finite
supports (also not known in advance). Second, another paper at this conference \citep{GaCa11} also deals with the
$\cK$--strategy which we study in Theorem~\ref{th:mainBer}; they however do not obtain second-order terms in closed forms
as we do and later extend their strategy to exponential families of distributions (while we extend our strategy
to the case of distributions with finite supports). On the other hand, they show how the $\cK$--strategy
can be extended in a straightforward manner to guarantee bounds with respect to the family of all bounded distributions
on a known interval; these bounds are suboptimal but improve on the ones of UCB-type algorithms.

\section{Definitions and tools}
\label{sec:Sanov}

Let $\cX$ be a Polish space; in the next sections, we will consider $\cX = \{ 0,1 \}$ or $\cX = [0,1]$.
We denote by $\cP(\cX)$ the set of probability distributions over $\cX$ and equip $\cP(\cX)$ with the distance $d$ induced by the norm $\norm$
defined by $\norm[\nu] = \sup_{f \in \cL} \, \bigl| \int_{\cX} f \,\d\nu \bigr|,$ where $\cL$ is the set of Lipschitz functions over $\cX$, taking values in $[-1,1]$
and with Lipschitz constant smaller than 1.

\paragraph{Kullback-Leibler divergence:} For two elements $\nu,\,\nup \in \cP(\cX)$, we write $\nu \ll \nup$ when $\nu$ is absolutely continuous
with respect to $\nup$ and denote in this case by $\d\nu/\d\nup$ the density of $\nu$ with respect to $\nup$.
We recall that the Kullback-Leibler divergence between $\nu$ and $\nup$ is defined
as
\beq\label{eq:def:KL}
\cK(\nu,\nup) = \int_{[0,1]} \frac{\d\nu}{\d\nup} \log \frac{\d\nu}{\d\nup} \, \d\nup \quad \mbox{if} \ \nu \ll \nup;
\qquad \mbox{and} \quad \cK(\nu,\nup) = +\infty \quad \mbox{otherwise.}
\eeq

\paragraph{Empirical distribution:} We consider a sequence $X_1,X_2,\ldots$ of random variables taking values in $\cX$, independent and identically distributed according to a distribution $\nu$. For all integers $t \geq 1$, we denote
the empirical distribution corresponding to the first $t$ elements of the sequence by
\[
\hat{\nu}_t = \frac{1}{t} \sum_{s=1}^t \delta_{X_t}\,.
\]

\paragraph{Non-asymptotic Sanov's Lemma:} The following lemma follows from a straightforward adaptation of~\citet[Theorem~2.1 and comments on page~372]{Din92}. \ifarXiv{Details of the proof are provided in the appendix.}\ifCOLT{Details of the proof are provided in
the extended version~\cite{Ext} of the present paper.}

\begin{lemma}
\label{lm:Sanov2}
Let $\cC$ be an open convex subset of $\cP(\cX)$ such that
\quad $
\displaystyle{\Lambda(\cC) = \inf_{\nup \in \cC} \, \cK(\nup,\nu) < \infty\,.}
$ \\
Then, for all $t \geq 1$, one has $\qquad \qquad \displaystyle{\P_\nu \bigl\{ \hat{\nu}_t \in \ol{\cC} \bigr\} \leq e^{- t \Lambda(\ol{\cC})}}$
\qquad where $\ol{\cC}$ is the closure of $\cC$.
\end{lemma}

This lemma should be thought of as a deviation inequality.
The empirical distribution converges (in distribution) to $\nu$. Now,
if (and only if) $\nu$ is not in the closure of $\cC$, then
$\Lambda(\cC) > 0$ and the lemma indicates how unlikely it is that
$\hat{\nu}_t$ is in this set $\ol{\cC}$ not containing the limit $\nu$.
The probability of interest decreases at a geometric rate, which depends on $\Lambda(\cC)$.

\section{Finite-time analysis for Bernoulli distributions}

In this section, we start with the case of Bernoulli distributions. Although this case is a special case of the general results of Section~\ref{sec:genberal.case}, we provide here a
complete and self-contained analysis of this case, where, in addition, we are able to provide closed forms for all the terms in the regret bound.
Note however that the resulting bound is slightly worse than what could be derived from the general case (for which more sophisticated tools are used).
This result is mainly provided as a warm-up.

\subsection{Reminder of some useful results for Bernoulli distributions}
\label{sec:reminderBer}

We denote by $\cB$ the subset of $\cP \bigl( [0,1] \bigr)$ formed by the Bernoulli
distributions; it corresponds to $\cB = \cP \bigl( \{0,1\} \bigr)$.
A generic element of $\cB$ will be denoted by $\beta(p)$, where
$p \in [0,1]$ is the probability mass put on $1$. We consider a sequence $X_1,X_2,\ldots$ of independent
and identically distributed random variables, with common distribution $\beta(p)$;
for the sake of clarity we will index, in this subsection only,
all probabilities and expectations with $p$.

For all integers $t \geq 1$, we denote by $\quad \displaystyle{\wh{p}_t = \frac{1}{t} \sum_{s=1}^t X_t} \quad$ the
empirical average of the first $t$ elements of the sequence.

The lemma below follows from an adaptation of~\citet[Proposition~2]{GaLe10}.
\ifarXiv{The details of the adaptation (and simplification) can be found in the appendix.}

\begin{lemma}
\label{lm:Aure}
For all $p \in [0,1]$, all $\eps > 1$, and all $t \geq 1$,
\[
\P_p \! \left( \bigcup_{s = 1}^t  \biggl\{ s\,\,\cK\Bigl( \beta\bigl({\wh{p}_s}\bigr), \, \beta(p) \Bigr) \geq \eps
\biggr\} \right) \leq 2e \, \bigl\lceil \eps \log t \bigr\rceil \, e^{-\eps}\,.
\]
In particular, for all random variables $N_t$ taking values in $\{ 1,\ldots,t \}$,
\[
\P_p \biggl\{ N_t \,\, \cK\Bigl( \beta\bigl({\wh{p}_{N_t}}\bigr), \, \beta(p) \Bigr) \geq \eps
\biggr\} \leq 2e \, \bigl\lceil \eps \log t \bigr\rceil \, e^{-\eps}\,.
\]
\end{lemma}

Another immediate fact about Bernoulli distributions is that for all $p \in (0,1)$,
the mappings $\cK_{\,\cdot\,,p} : q \in (0,1) \mapsto \cK \bigl( \beta(p),\beta(q) \bigr)$ and
$\cK_{p,\,\cdot\,} : q \in [0,1] \mapsto \cK \bigl( \beta(q), \beta(p) \bigr)$ are continuous and take finite values.
In particular, we have, for instance, that for all $\eps > 0$ and $p \in (0,1)$, the set
\[
\Bigl\{ q \in [0,1] : \ \ \cK \bigl( \beta(p),\beta(q) \bigr) \leq \eps \Bigr\}
\]
is a closed interval containing $p$. This property still holds when $p \in \{ 0,1 \}$,
as in this case, the interval is reduced to $\{p\}$.

\subsection{Strategy and analysis}

We consider the so-called {\em $\cK$--strategy} of Figure~\ref{fig:Ber}, which was already considered
in the literature, see \cite{Burnetas96,Filipi2010}.
The numerical computation of the quantities $B^+_{a,t}$ is straightforward (by
convexity of $\cK$ in its second argument, by using iterative methods) and is detailed therein.

\begin{figure}[t]
\rule{\linewidth}{.5pt}
{\small
\emph{Parameters}: A non-decreasing function $f : \mathbb{N} \to \mathbb{R}$ \medskip \\
\emph{Initialization}: Pull each arm of $\cA$ once \medskip \\
\emph{For} rounds $t+1$, where $t \geq |\cA|$,
\begin{itemize}
\item[--] compute for each arm $a \in \cA$ the quantity
\[
B^+_{a,t} = \max \, \biggl\{ q \in [0,1] : \ \ N_t(a) \,\, \cK\Bigl( \beta\bigl( \hat{\mu}_{a,N_t(a)} \bigr), \, \beta(q) \Bigr) \leq
f(t) \biggr\}\,,
\]
where $\qquad \qquad \displaystyle{ \hat{\mu}_{a,N_t(a)} = \frac{1}{N_t(a)} \sum_{s \leq t: \, A_s = a} Y_s\,; }$
\item[--] in case of a tie, pick an arm with largest value of $\hat{\mu}_{a,N_t(a)}$;
\item[--] pull any arm $\displaystyle{ A_{t+1} \in \mathop{\mathrm{argmax}}_{a \in \cA} \, B^+_{a,t}\,. }$ \vspace{-.25cm}
\end{itemize}
}
\rule{\linewidth}{.5pt} \vspace{-.5cm}
\caption{\label{fig:Ber} The $\cK$--strategy.}
\end{figure}

Before proceeding, we denote by $\sigma^2_a = \mu_a (1-\mu_a)$ the variance
of each arm $a \in \cA$ (and take the short-hand notation $\sigma^{\star,2}$
for the variance of an optimal arm).

\begin{theorem}
\label{th:mainBer}
When $\mu^\star \in (0,1)$,
for all non-decreasing functions $f : \mathbb{N} \to \mathbb{R}_+$ such that $f(1) \geq 1$,
the expected regret
$R_T$ of the strategy of Figure~\ref{fig:Ber} is upper bounded by the infimum, as the $(c_a)_{a\in\cA}$ describe $(0,+\infty)$,
of the quantities
\[
\sum_{a \in \cA} \Delta_a
\Biggl( \frac{(1+c_a)\,f(T)}{\cK\bigl( \beta(\mu_a), \, \beta(\mu^\star) \bigr)}
+ 4 e \sum_{t = |\cA|}^{T-1} \bigl\lceil f(t) \log t \bigr\rceil \, e^{-f(t)}
+ \frac{(1+c_a)^2}{8 \, c_a^2
\Delta_a^2 \, \min \bigl\{ \sigma_a^4,\,\sigma^{\star,4} \bigr\} } \mathbb{I}_{ \{ \mu_a \in (0,1) \} }
+3 \Biggr)\,.
\]
For $\mu^\star = 0$, its regret is null. For $\mu^\star = 1$, it satisfies $R_T \leq 2 \bigl( |\cA| - 1 \bigr)$.
\end{theorem}

A possible choice for the function $f$ is $f(t) = \log \bigl( (et) \log^3(et) \bigr)$, which is non decreasing,
satisfies $f(1) \geq 1$, and is such that the second term in the sum above is bounded (by a basic result about
so-called Bertrand's series). Now, as
the constants $c_a$ in the bound are parameters of the analysis (and not of the strategy), they can be optimized.
For instance, with the choice of $f(t)$ mentioned above,
taking each $c_a$ proportional to $(\log T)^{-1/3}$ (up to a multiplicative constant that depends on
the distributions $\nu_a$) entails the regret bound
\[
\sum_{a \in \cA} \Delta_a
\frac{\log T}{\cK\bigl( \beta(\mu_a), \, \beta(\mu^\star) \bigr)}
+ \varepsilon_T\,,
\]
where it is easy to give an explicit and closed-form expression of $\varepsilon_T$; in this
conference version, we only indicate that $\varepsilon_T$ is of order of $(\log T)^{2/3}$
but we do not know whether the order of magnitude of this second-order term is optimal. \\

\begin{proof}
We first deal with the case where $\mu^{\star} \not\in \{ 0,1 \}$ and introduce an additional notation. In view of the remark at the end of Section~\ref{sec:reminderBer},
for all arms $a$ and rounds $t$, we let $B^-_{a,t}$ be the element in $[0,1]$ such that
\begin{equation}
\label{eq:defK}
\biggl\{ q \in [0,1] : \ \ N_t(a) \,\, \cK\Bigl( \beta\bigl( \hat{\mu}_{a,N_t(a)} \bigr), \, \beta(q) \Bigr) \leq
f(t) \biggr\} = \bigl[ B^-_{a,t}, \,\, B^+_{a,t} \bigr]\,.
\end{equation}
As~(\ref{eq:defregr}) indicates, it suffices to bound $N_T(a)$ for all suboptimal arms $a$, i.e.,
for all arms such that $\mu_a < \mu^\star$. We will assume in addition that $\mu_a > 0$ (and we also
have $\mu_a \leq \mu^\star < 1$); the case where $\mu_a = 0$ will be handled separately. \smallskip \\
\indent \textbf{Step 1: A decomposition of the events of interest.}
For $t \geq |\cA|$, when $A_{t+1} = a$, we have in particular,
by definition of the strategy, that $B^+_{a,t} \geq B^+_{a^\star,t}$.
On the event
\[
\bigl\{ A_{t+1} = a \bigr\} \, \cap \,
\Bigl\{ \mu^\star \in \bigl[ B^-_{a^\star,t}, \,\, B^+_{a^\star,t} \bigr] \Bigr\} \, \cap \,
\Bigl\{ \mu_a \in \bigl[ B^-_{a,t}, \,\, B^+_{a,t} \bigr] \Bigr\}\,,
\]
we therefore have, on the one hand, $\mu^\star \leq B^+_{a^\star,t} \leq B^+_{a,t}$
and on the other hand, $B^-_{a,t} \leq \mu_a \leq \mu^{\star}$, that is, the considered
event is included in $\Bigl\{ \mu^\star \in \bigl[ B^-_{a,t}, \,\, B^+_{a,t} \bigr] \Bigr\}$.
We thus proved that
\[
\bigl\{ A_{t+1} = a \bigr\} \subseteq
\Bigl\{ \mu^\star \not\in \bigl[ B^-_{a^\star,t}, \,\, B^+_{a^\star,t} \bigr] \Bigr\} \, \cup \,
\Bigl\{ \mu_a \not\in \bigl[ B^-_{a,t}, \,\, B^+_{a,t} \bigr] \Bigr\} \, \cup \,
\Bigl\{ \mu^\star \in \bigl[ B^-_{a,t}, \,\, B^+_{a,t} \bigr] \Bigr\}\,.
\]
Going back to the definition~(\ref{eq:defK}), we get in particular the inclusion
\begin{align*}
\bigl\{ A_{t+1} = a \bigr\} \subseteq & \quad
\biggl\{ N_t(a^\star) \,\, \cK\Bigl( \beta\bigl( \hat{\mu}_{a^\star,N_t(a^\star)} \bigr), \, \beta(\mu^\star) \Bigr) >
f(t) \biggr\} \\
& \cup \, \biggl\{ N_t(a) \,\, \cK\Bigl( \beta\bigl( \hat{\mu}_{a,N_t(a)} \bigr), \, \beta(\mu_a) \Bigr) >
f(t) \biggr\} \\
& \cup \, \Biggl( \biggl\{ N_t(a) \,\, \cK\Bigl( \beta\bigl( \hat{\mu}_{a,N_t(a)} \bigr), \, \beta(\mu^\star) \Bigr) \leq
f(t) \biggr\} \, \cap \, \bigl\{ A_{t+1} = a \bigr\} \Biggr)\,. \smallskip
\end{align*}

\textbf{Step 2: Bounding the probabilities of two elements of the decomposition.}
We consider the filtration $(\cF_t)$, where for all $t \geq 1$, the $\sigma$--algebra $\cF_t$ is generated
by $A_1,Y_1$,\,$\ldots$,\, $A_t,Y_t$. In particular, $A_{t+1}$ and thus all $N_{t+1}(a)$ are $\cF_t$--measurable.
We denote by $\tau_{a,1}$ the deterministic round at which $a$ was pulled for the first time
and by $\tau_{a,2},\,\tau_{a,3},\,\ldots$ the rounds $t \geq |\cA|+1$ at which $a$ was then played; since
for all $k \geq 2$,
\[
\tau_{a,k} = \min \bigl\{ t \geq |\cA|+1 : \ \ N_t(a) = k \bigr\}\,,
\]
we see that $\bigl\{ \tau_{a,k} = t \bigr\}$ is $\cF_{t-1}$--measurable. Therefore, for each $k \geq 1$, the
random variable $\tau_{a,k}$ is a (predictable) stopping time. Hence, by a well-known fact in probability theory
(see, e.g., \citealt[Section~5.3]{ChTe88}), the random variables $\widetilde{X}_{a,k} = Y_{\tau_{a,k}}$, where $k = 1,2,\ldots$
are independent and identically distributed according to $\nu_a$.
Since on $\bigl\{ N_t(a) = k \bigr\}$, we have the rewriting
\[
\hat{\mu}_{a,N_t(a)} = \widetilde{\mu}_{a,k}\, \qquad \mbox{where} \qquad \widetilde{\mu}_{a,k} = \frac{1}{k} \sum_{j=1}^k \widetilde{X}_{a,j}\,,
\]
and since for $t \geq |\cA|+1$, one has $N_t(a) \geq 1$ with probability 1,
we can apply the second statement in Lemma~\ref{lm:Aure} and get, for all $t \geq |\cA|+1$,
\[
\P \biggl\{ N_t(a) \,\, \cK\Bigl( \beta\bigl( \hat{\mu}_{a,N_t(a)} \bigr), \, \beta(\mu_a) \Bigr) >
f(t) \biggr\} \leq 2e \, \bigl\lceil f(t) \log t \bigr\rceil \, e^{-f(t)}\,.
\]
A similar argument shows that for all $t \geq |\cA|+1$,
\[
\P \biggl\{ N_t(a^\star) \,\, \cK\Bigl( \beta\bigl( \hat{\mu}_{a^\star,N_t(a^\star)} \bigr), \, \beta(\mu^\star) \Bigr) >
f(t) \biggr\} \leq 2e \, \bigl\lceil f(t) \log t \bigr\rceil \, e^{-f(t)}\,. \smallskip
\]

\textbf{Step 3: Rewriting the remaining terms.}
We therefore proved that
{\small
\begin{align*}
\E \bigl[ N_T(a) \bigr]
\leq 1 & + 4 e \sum_{t = |\cA|}^{T-1} \bigl\lceil f(t) \log t \bigr\rceil \, e^{-f(t)}
+ \sum_{t = |\cA|}^{T-1} \P \Biggl( \biggl\{ N_t(a) \,\, \cK\Bigl( \beta\bigl( \hat{\mu}_{a,N_t(a)} \bigr), \, \beta(\mu^\star) \Bigr) \leq
f(t) \biggr\} \, \cap \, \bigl\{ A_{t+1} = a \bigr\} \Biggr)
\end{align*}
}
and deal now with the last sum. Since $f$ is non decreasing, it is bounded by
\[
\sum_{t = |\cA|}^{T-1} \, \P \Bigl( K_t \, \cap \, \bigl\{ A_{t+1} = a \bigr\} \Bigr) \qquad
\mbox{where} \qquad
K_t = \biggl\{ N_t(a) \,\, \cK\Bigl( \beta\bigl( \hat{\mu}_{a,N_t(a)} \bigr), \, \beta(\mu^\star) \Bigr) \leq
f(T) \biggr\} \, .
\]
Now, $\qquad \displaystyle{ \sum_{t = |\cA|}^{T-1} \, \P \Bigl( K_t \, \cap \, \bigl\{ A_{t+1} = a \bigr\} \Bigr)
= \E \! \left[ \sum_{t = |\cA|}^{T-1} \ind_{ \bigl\{ A_{t+1} = a \bigr\} } \ind_{K_t} \right]
= \E \! \left[ \sum_{k \geq 2} \ind_{ \bigl\{ \tau_{a,k} \leq T \bigr\} } \ind_{K_{\tau_{a,k}-1}} \right].
}$ \vspace{.15cm} \\
We note that, since $N_{\tau_{a,k} - 1}(a) = k-1$, we have that
\[
K_{\tau_{a,k}-1} =
\biggl\{ (k-1) \,\, \cK\Bigl( \beta\bigl( \widetilde{\mu}_{a,k-1} \bigr), \, \beta(\mu^\star) \Bigr) \leq f(T) \biggr\}\,.
\]
All in all, since $\tau_{a,k} \leq T$ implies $k \leq T-|\cA|+1$ (as each arm is played at least once during the
first $|\cA|$ rounds), we have
{\small
\begin{equation}
\label{eq:step3} \!\!
\E \! \left[ \sum_{k \geq 2} \ind_{ \bigl\{ \tau_{a,k} \leq T \bigr\} } \ind_{K_{\tau_{a,k}-1}} \right]
\leq \E \! \left[ \sum_{k = 2}^{T-|\cA|+1} \ind_{K_{\tau_{a,k}-1}} \right]
= \! \sum_{k = 2}^{T-|\cA|+1} \P \biggl\{ (k-1) \,\, \cK\Bigl( \beta\bigl( \widetilde{\mu}_{a,k-1} \bigr), \, \beta(\mu^\star) \Bigr) \leq f(T) \biggr\}\,. \smallskip
\end{equation}
}

\textbf{Step 4: Bounding the probabilities of the latter sum via Sanov's lemma.}
For each $\gamma > 0$, we define the convex open set
$\displaystyle{\cC_\gamma = \Bigl\{ \beta(q) \in \cB : \ \ \cK\bigl( \beta(q), \, \beta(\mu^\star) \bigr) < \gamma \Bigr\}
},$
which is a non-empty set (since $\mu^\star < 1$);
by continuity of the mapping $\cK_{\,\cdot\,,\mu^\star}$ defined after
the statement of Lemma~\ref{lm:Aure} when $\mu^\star \in (0,1)$,
its closure equals
$\displaystyle{
\ol{\cC}_\gamma = \Bigl\{ \beta(q) \in \cB : \ \ \cK\bigl( \beta(q), \, \beta(\mu^\star) \bigr) \leq \gamma \Bigr\}\,.
}$

In addition, since $\mu_a \in (0,1)$, we have that $\cK\bigl( \beta(q), \, \beta(\mu_a) \bigr) < \infty$ for all
$q \in [0,1]$. In particular, for all $\gamma > 0$,
the condition $\Lambda \bigl( \cC_\gamma \bigr) < \infty$ of Lemma~\ref{lm:Sanov2} is satisfied.
Denoting this value by
\[
\theta_a(\gamma) = \inf \biggl\{ \cK \bigl( \beta(q), \, \beta(\mu_a) \bigr) : \ \
\beta(q) \in \cB \ \ \mbox{\small such that} \ \ \cK\bigl( \beta(q), \, \beta(\mu^\star) \bigr) \leq \gamma \biggr\}\,,
\]
we get by the indicated lemma that for all $k \geq 1$,
\[
\P \biggl\{ \cK\Bigl( \beta\bigl( \widetilde{\mu}_{a,k} \bigr), \, \beta(\mu^\star) \Bigr) \leq \gamma \biggr\}
= \P \Bigl\{ \beta\bigl( \widetilde{\mu}_{a,k} \bigr) \in \ol{\cC}_\gamma \Bigr\}
\leq e^{-k \, \theta_a(\gamma)}\,.
\]
Now, since (an open neighborhood of) $\beta(\mu_a)$ is not included in $\ol{\cC}_\gamma$ as soon as
$0 < \gamma < \cK\bigl( \beta(\mu_a), \, \beta(\mu^\star) \bigr)$, we have that
$\theta_a(\gamma) > 0$ for such values of $\gamma$.
To apply the obtained inequality to the last sum in~(\ref{eq:step3}), we fix a constant $c_a > 0$ and denote by $k_0$
the following upper integer part,
$\displaystyle{
k_0 = \left\lceil \frac{(1+c_a)\,f(T)}{\cK\bigl( \beta(\mu_a), \, \beta(\mu^\star) \bigr)} \right\rceil,
}$
so that $f(T)/k \leq \cK\bigl( \beta(\mu_a), \, \beta(\mu^\star) \bigr)/(1+c_a) <
\cK\bigl( \beta(\mu_a), \, \beta(\mu^\star) \bigr)$ for $k \geq k_0$, hence,
\begin{eqnarray*}
\sum_{k = 2}^{T-|\cA|+1} \, \P \biggl\{ (k-1) \,\, \cK\Bigl( \beta\bigl( \widetilde{\mu}_{a,k-1} \bigr), \, \beta(\mu^\star) \Bigr) \leq f(T) \biggr\}
& \leq & \sum_{k = 1}^{T} \,
\P \biggl\{ \cK\Bigl( \beta\bigl( \widetilde{\mu}_{a,k} \bigr), \, \beta(\mu^\star) \Bigr) \leq \frac{f(T)}{k} \biggr\} \\
& \leq & k_0-1 + \sum_{k = k_0}^{T} \, \exp \Bigl( -k \, \theta_a\bigl( f(T)/k \bigr) \Bigr)\,.
\end{eqnarray*}
Since $\theta_a$ is a non-increasing function,
\begin{eqnarray*}
\sum_{k = k_0}^{T} \, \exp \Bigl( -k \, \theta_a\bigl( f(T)/k \bigr) \Bigr)
& \leq & \sum_{k = k_0}^{T} \, \exp \Bigl( -k \, \theta_a\bigl( \cK\bigl( \beta(\mu_a), \, \beta(\mu^\star) \bigr)/(1+c_a) \bigr) \Bigr) \\
& \leq & \Gamma_a(c_a) \, \exp \Bigl( -k_0 \, \theta_a\bigl( \cK\bigl( \beta(\mu_a), \, \beta(\mu^\star) \bigr)/(1+c_a) \bigr) \Bigr) \leq \Gamma_a(c_a),
\end{eqnarray*}
where $\displaystyle{
\Gamma_a(c_a) = \Big[ 1 - \exp \Bigl( - \theta_a\bigl( \cK\bigl( \beta(\mu_a), \, \beta(\mu^\star) \bigr)/(1+c_a) \bigr) \Bigr) \Big]^{-1}}\,.
$ \\
Putting all pieces together, we thus proved so far that
\[
\E \bigl[ N_T(a) \bigr]
\leq 1 + \frac{(1+c_a)\,f(T)}{\cK\bigl( \beta(\mu_a), \, \beta(\mu^\star) \bigr)}
+ 4 e \sum_{t = |\cA|}^{T-1} \bigl\lceil f(t) \log t \bigr\rceil \, e^{-f(t)} + \Gamma_a(c_a)
\]
and it only remains to deal with $\Gamma_a(c_a)$. \medskip \\
\indent \textbf{Step 5: Getting an upper bound in closed form for $\Gamma_a(c_a)$.}
We will make repeated uses of Pinsker's inequality: for $p,q \in [0,1]$,
one has $\cK \bigl( \beta(p),\beta(q) \bigr) \geq 2 \, (p-q)^2\,.$ \\
In what follows, we use the short-hand notation
$\Theta_a = \theta_a\bigl( \cK\bigl( \beta(\mu_a), \, \beta(\mu^\star) \bigr)/(1+c_a) \bigr)$
and therefore need to upper bound $1/\bigl( 1 - e^{- \Theta_a} \bigr)$.
Since for all $u \geq 0$, one has $e^{-u} \leq 1 - u +u^2/2$, we get
$\displaystyle{
\Gamma_a(c_a) \leq \frac{1}{\Theta_a \bigl( 1 - \Theta_a/2 \bigr)} \leq \frac{2}{\Theta_a}
 \mbox{ for }  \Theta_a \leq 1,
}$
and
$\displaystyle{
\Gamma_a(c_a) \leq \frac{1}{1-e^{-1}} \leq 2  \mbox{ for }  \Theta_a \geq 1.
}$
It thus only remains to lower bound $\Theta_a$ in the case when it is smaller than 1.

By the continuity properties of the Kullback-Leibler divergence, the infimum in the definition
of $\theta_a$ is always achieved; we therefore let $\mt$ be an element in $[0,1]$
such that
\[
\Theta_a = \cK \bigl( \beta({\widetilde{\mu}}), \, \beta({\mu_a}) \bigr)
\qquad \mbox{and} \qquad
\cK \bigl( \beta({\widetilde{\mu}}), \, \beta({\mu^{\star}}) \bigr) =
\frac{\cK \bigl( \beta({\mu_a}), \, \beta({\mu^{\star}}) \bigr)}{1+c}\,;
\]
it is easy to see that we have the ordering $\mu_a < \widetilde{\mu} < \mu^{\star}$.
By Pinsker's inequality, $\Theta_a \geq 2 \bigl( \widetilde{\mu} - \mu_a \bigr)^2$
and we now lower bound the latter quantity. We use the short-hand notation
$f(p) = \cK\bigl(\beta(p),\beta({\mu^\star})\bigr)$ and note that the thus defined mapping
$f$ is convex and differentiable on $(0,1)$;
its derivative equals
$
f'(p) = \log \bigl( (1-\mu^{\star}) / (\mu^\star) \bigr) + \log \bigl( p/(1-p) \bigr)
$
for all $p \in (0,1)$ and is therefore non positive for $p \leq \mu^\star$.
By the indicated convexity of $f$, using a sub-gradient inequality, we get
$
f\bigl(\mt\bigr) - f(\mu_a) \geq f'(\mu_a) \, \bigl( \mt - \mu_a \bigr)\,,
$
which entails, since $f'(\mu_a) < 0$,
\begin{equation}
\label{eq:Theta}
\mt - \mu_a \geq \frac{f\bigl(\mt\bigr) - f(\mu_a)}{f'(\mu_a)} = \frac{c_a}{1+c_a} \,\, \frac{f(\mu_a)}{-f'(\mu_a)}\,,
\end{equation}
where the equality follows from the fact that by definition of $\mu$, we have $f\bigl(\mt\bigr) = f(\mu_a)/(1+c_a)$.
Now, since $f'$ is differentiable as well on $(0,1)$ and takes the value $0$ at $\mu^\star$, a Taylor's equality entails that there exists
a $\xi \in (\mu_a,\mu^\star)$ such that
\[
-f'(\mu_a) = f'(\mu^\star) - f'(\mu_a) = f''(\xi) \, \bigl( \mu^\star - \mu_a)
\qquad \mbox{where} \quad f''(\xi) = 1/\xi + 1/(1-\xi) = 1 \big/ \bigl( \xi(1-\xi) \bigr)\,.
\]
Therefore, by convexity of $\tau \mapsto \tau(1-\tau)$, we get that
\[
\frac{1}{-f'(\mu_a)} \geq \frac{\min \bigl\{ \mu_a(1-\mu_a),\,\mu^{\star}(1-\mu^\star) \bigr\}}{\mu^\star - \mu_a}\,.
\]
Substituting this into~(\ref{eq:Theta}) and using again Pinsker's inequality to lower bound $f(\mu_a)$, we
have proved
\[
\mt - \mu_a \geq 2 \, \frac{c_a}{1+c_a} \,
\bigl( \mu^\star - \mu_a \bigr) \, \min \bigl\{ \mu_a(1-\mu_a),\,\mu^{\star}(1-\mu^\star) \bigr\}\,.
\]
Putting all pieces together, we thus proved that
\[
\Gamma_a(c_a) \leq 2 \, \max \left\{ \frac{(1+c_a)^2}{8 \, c_a^2
\bigl( \mu^\star - \mu_a \bigr)^2 \, \Bigl( \min \bigl\{ \mu_a(1-\mu_a),\,\mu^{\star}(1-\mu^\star) \bigr\} \Bigr)^2 }
, \,\,1 \right\}\,;
\]
bounding the maximum of the two quantities by their sum concludes the main part of the proof.
\smallskip \\

\textbf{Step 6: For $\mu^{\star} \in \{ 0,1 \}$ and/or $\mu_a = 0$.}
When $\mu^{\star} = 1$, then $\hat{\mu}_{a^\star,N_t(a\star)} = 1$ for all $t \geq |\cA|+1$,
so that $B^+_{a^\star,t} = 1$ for all $t \geq |\cA|+1$.
Thus, the arm $a$ is played after round $t \geq |\cA|+1$
only if $B^+_{a,t} = 1$ and $\hat{\mu}_{a,N_t(a)} = 1$ (in view of the
tie-breaking rule of the considered strategy). But this means that
$a$ is played as long as it gets payoffs equal to 1 and is stopped being played
when it receives the payoff 0 for the first time. Hence, in this case, we have that
the sum of payoffs equals at least $T-2\bigl( |\cA| - 1)$ and the regret $R_T = \E[ T \mu^\star - (Y_1+\ldots+Y_t)]$ is therefore
bounded by $2\bigl( |\cA| - 1)$.

When $\mu^\star = 0$, a Dirac mass over 0 is associated with all arms and the regret
of all strategies is equal to 0.

We consider now the case $\mu^{\star} \in (0,1)$ and $\mu_a = 0$, for which the first three steps
go through; only in the upper bound of step 4 we used the fact that $\mu_a > 0$.
But in this case, we have a deterministic bound on~(\ref{eq:step3}).
Indeed, since $\cK \bigl( \beta(0), \beta(\mu^\star) \bigr) = - \log \mu^\star$, we have
$k \, \cK \bigl( \beta(0), \beta(\mu^\star) \bigr) \leq f(T)$ if and only if
\[
k \leq \frac{f(T)}{- \log \mu^\star} = \frac{f(T)}{\cK \bigl( \beta(\mu_a), \beta(\mu^\star) \bigr)}\,,
\]
which improves on the general bound exhibited in step 4.
\end{proof}

\begin{remark} {\em
Note that Step~5 in the proof is specifically designed to provide
an upper bound on $\Gamma_a(c_a)$ in the case of Bernoulli distributions.
In the general case, getting such an explicit bound seems more involved. }
\end{remark}

\section{A finite-time analysis in the case of distributions with finite support}
\label{sec:genberal.case}
Before stating and proving our main result, Theorem~\ref{thm:main}, we introduce the quantity $\kmin$ and list some of its properties.

\subsection{Some useful properties of $\kmin$ and its level sets}
\label{sec:setkmin}

We now introduce the key quantity in order to generalize the previous algorithm to handle the case of distributions with finite support.
To that end, we introduce $\cP_F \bigl( [0,1] \bigr)$,
the subset of $\cP \bigl( [0,1] \bigr)$ that consists of distributions with finite support.

\begin{definition}
For all distributions $\nu \in \cP_F \bigl( [0,1] \bigr)$ and $\mu \in [0,1)$, we define
\beqan
\kmin(\nu,\mu) = \inf \, \Bigl\{ \cK(\nu,\nu') : \ \ \nu' \in \cP_F \bigl( [0,1] \bigr) \ \ \mbox{\rm s.t.} \ \ E(\nu') > \mu \Bigr\},
\eeqan
where $E(\nu')= \int_{[0,1]} x \, {\mbox{\rm d}}\nu'(x)$ denotes the expectation of the distribution $\nu'$.
\end{definition}

%

We now remind some useful properties of $\kmin$. \citet[Lemma~6]{HondaArxiv}
can be reformulated in our context as follows.

\begin{lemma}\label{lem:sci}
For all $\nu \in \cP_F \bigl( [0,1] \bigr)$, the  mapping $\kmin(\nu,\,\cdot\,)$ is continuous and non decreasing in its argument $\mu \in [0,1)$.
Moreover, the mapping $\kmin(\,\cdot\,,\mu)$ is lower semi-continuous on $\cP_F \bigl( [0,1] \bigr)$ for all $\mu \in [0,1)$.
\end{lemma}

The next two lemmas bound the variation of $\kmin$, respectively in its first and second arguments.
(For clarity, we denote the expectations with respect to $\nu$ by $\E_\nu$.)
Their proofs \ifCOLT{can be found in the extended version of the present conference paper~\citep{Ext}.}\ifarXiv{are both deferred to the appendix.}
We denote by $\norm_1$ the $\ell^1$--norm on $\cP \bigl( [0,1] \bigr)$ and recall that
the $\ell^1$--norm of $\nu-\nu'$ corresponds to twice the distance in variation between $\nu$ and $\nu'$.

\begin{lemma}
\label{lm:pinsker}
For all $\mu \in (0,1)$ and for all $\nu,\,\nu'\in \cP_F \bigl( [0,1] \bigr)$, the following holds true.
\begin{itemize}
\item[--] In the case when $\E_\nu \bigl[ (1-\mu)/(1-X) \bigr] > 1$, then
$
\kmin(\nu,\mu) - \kmin(\nu',\mu) \leq M_{\nu,\mu} \, \Arrowvert \nu - \nu' \Arrowvert_1\,,
$
for some constant $M_{\nu,\mu} > 0$.
\item[--] In the case when $\E_\nu \bigl[ (1-\mu)/(1-X) \bigr] \leq 1$,
the fact that $\kmin(\nu,\mu) - \kmin(\nu',\mu) \geq \alpha \, \kmin(\nu,\mu)$
for some $\alpha \in (0,1)$ entails that
\[
\Arrowvert \nu - \nu' \Arrowvert_1 \geq
\frac{1-\mu}{(2/\alpha) \, \bigl( (2/\alpha) - 1 \bigr)}\,.
\]

\end{itemize}
\end{lemma}

\begin{lemma}
\label{lm:dminvar}
We have that
for any $\nu \in \cP_F \bigl( [0,1] \bigr)$, provided that $\mu\geq \mu-\epsilon > E(\nu)$, the following inequalities hold true:
\beqan
\epsilon/ (1-\mu) \geq \kmin( \nu,\mu) - \kmin( \nu,\mu - \epsilon) \geq  2\epsilon^2
\eeqan
Moreover, the first inequality is also valid when $E(\nu) \geq \mu > \mu-\epsilon$
or $\mu > E(\nu) \geq \mu-\epsilon$.
\end{lemma}

\paragraph{Level sets of $\kmin$:} For each $\gamma > 0$ and $\mu \in (0,1)$, we consider the set
\begin{eqnarray*}
\cC_{\mu,\gamma} & = & \Bigl\{ \nu' \in \cP_F \bigl( [0,1] \bigr) : \ \ \kmin( \nu',\mu) < \gamma \Bigr\} \\
& = & \Bigl\{ \nu' \in \cP_F \bigl( [0,1] \bigr) : \ \ \exists \, \nu'_\mu \in \cP_F \bigl( [0,1] \bigr)
\ \ \mbox{s.t.} \ \ E\bigl(\nu'_\mu\bigr) > \mu \ \ \mbox{and} \ \ \cK\bigl( \nu',\nu'_\mu\bigr) < \gamma \Bigr\}\,.
\end{eqnarray*}
We detail a property in the following lemma, \ifarXiv{whose proof is also deferred to the appendix.}\ifCOLT{whose proof
can be found in the extended version of the present conference paper~\citep{Ext}.}

\begin{lemma}
\label{lm:closure}
For all $\gamma > 0$ and $\mu \in (0,1)$, the
set $\cC_{\mu,\gamma}$ is a non-empty open convex set. Moreover,
\beqan
\ol{\cC}_{\mu,\gamma} \, \supseteq \, \Bigl\{ \nu' \in \cP_F \bigl( [0,1] \bigr) : \ \ \kmin( \nu',\mu) \leq \gamma \Bigr\}\,.
\eeqan
\end{lemma}

\subsection{The $\kmin$--strategy and a general performance guarantee}

For each arm $a \in \cA$ and round $t$ with $N_t(a) > 0$, we denote by $\hat{\nu}_{a,N_t(a)}$ the empirical
distribution of the payoffs obtained till round $t$ when picking arm $a$, that is,
\[
\hat{\nu}_{a,N_t(a)} = \frac{1}{N_t(a)} \sum_{s \leq t: \, A_s = a} \delta_{Y_s}\,,
\]
where for all $x \in [0,1]$, we denote by $\delta_x$ the Dirac mass on $x$.
We define the corresponding empirical averages as
\[
\hat{\mu}_{a^\star,N_t(a^\star)} = E \bigl( \hat{\nu}_{a^\star,N_t(a^\star)} \bigr) =
\frac{1}{N_t(a)} \sum_{s \leq t: \, A_s = a} Y_s\,.
\]
We then consider the {\em $\kmin$--strategy} defined in Figure~\ref{fig:Gal}. Note that the use of maxima in the definitions of the $B^+_{a,t}$ is justified by Lemma~\ref{lem:sci}.

As explained in~\cite{HondaArxiv}, the computation of the quantities $\kmin$ can be done efficiently in this case, i.e.,
when we consider only distributions with finite supports. This is because in the computation of $\kmin$, it is sufficient to consider only distributions with the same support as the empirical distributions (up to one point). Note that the knowledge of the support of the distributions associated with the arms is not required.

\begin{figure}[t]
\rule{\linewidth}{.5pt}
{\small
\emph{Parameters}: A non-decreasing function $f : \mathbb{N} \to \mathbb{R}$ \medskip \\
\emph{Initialization}: Pull each arm of $\cA$ once \medskip \\
\emph{For} rounds $t+1$, where $t \geq |\cA|$,
\begin{itemize}
\item[--] compute for each arm $a \in \cA$ the quantity
\[
B^+_{a,t} = \max \, \Bigl\{ q \in [0,1] : \ \ N_t(a) \,\, \kmin \bigl( \hat{\nu}_{a,N_t(a)}, \, q \bigr) \leq
f(t) \Bigr\}\,,
\]
where $\qquad \qquad \displaystyle{ \hat{\nu}_{a,N_t(a)} = \frac{1}{N_t(a)} \sum_{s \leq t: \, A_s = a} \delta_{Y_s}\,; }$ \vspace{-.25cm}
\item[--] in case of a tie, pick an arm with largest value of $\hat{\mu}_{a,N_t(a)}$;
\item[--] pull any arm $\displaystyle{ A_{t+1} \in \mathop{\mathrm{argmax}}_{a \in \cA} \, B^+_{a,t}\,. }$ \vspace{-.25cm}
\end{itemize}
}
\rule{\linewidth}{.5pt} \vspace{-.5cm}
\caption{\label{fig:Gal} The strategy $\kmin$.}
\end{figure}

\begin{theorem}
\label{thm:main}
Assume that $\nu^\star$ is finitely supported, with expectation $\mu^\star \in (0,1)$
and with support denoted by $\cS^\star$.
Let $a \in \cA$ be a suboptimal arm such that $\mu_a > 0$ and $\nu_a$ is finitely supported.
Then, for all $c_a > 0$ and all
\[
0 < \epsilon < \min \left\{ \Delta_a, \, \frac{c_a/2}{1+c_a} (1-\mu^\star) \, \kmin(\nu_a,\mu^\star) \right\},
\]
the expected number of times the $\kmin$--strategy, run with $f(t)=\log t$, pulls arm $a$ satisfies
\[
\E \bigl[ N_T(a) \bigr] \leq 1 +\frac{(1+c_a)\,\log T}{\kmin(\nu_a,\mu^\star)}
+ \frac{1}{1 - e^{- \Theta_a(c_a,\epsilon)}}
+ \frac{1}{\epsilon^2}\log\biggl(\frac{1}{1-\mu^* + \epsilon}\biggr) \sum_{k=1}^T (k+1)^{|\cS^\star|} \, e^{-k\epsilon^2}
+ \frac{1}{(\Delta_a-\epsilon)^2}\,,
\]
where
\[
\Theta_a(c_a,\epsilon) = \theta_a \! \left( \frac{\log T}{k_0} + \frac{\epsilon}{1-\mu^\star} \right)
\qquad \mbox{with} \qquad k_0 = \left\lceil \frac{(1+c_a)\,\log T}{\kmin(\nu_a,\mu^\star)} \right\rceil\,.
\]
and for all $\gamma > 0$,
\[
\theta_a(\gamma) = \inf \Bigl\{ \cK(\nu',\nu_a) : \ \ \nu' \ \,\, \mbox{\small \rm s.t.} \ \,\, \kmin(\nu',\mu^\star) < \gamma \Bigr\}\,.
\]
\end{theorem}

As a corollary, we get (by taking some common value for all $c_a$) that for all $c>0$,
\[
\ol{R}_T \leq \sum_{a\in\cA} \Delta_a \frac{(1+c)\,\log T}{\kmin(\nu_a,\mu^\star)} + h(c)\,,
\]
where $h(c)<\infty$ is a function of $c$ (and of the distributions associated with the arms),
which is however independent of $T$. As a consequence, we recover the asymptotic results of \citet{Burnetas96,Honda10}, i.e.,
the guarantee that
\[
\limsup_{T\rightarrow\infty} \frac{\ol{R}_T}{\log T} \leq \sum_{a\in\cA}  \frac{\Delta_a}{\kmin(\nu_a,\mu^\star)}\,.
\]

Of course, a sharper optimization can be performed
by carefully choosing the constants $c_a$, that are parameters of the analysis;
similarly to the comments after the statement of Theorem~\ref{th:mainBer}, we would then get
a dominant term with a constant factor $1$ instead of $1+c$ as above,
plus an additional second-order term. Details are left to a journal version of this paper.

\medskip
\begin{proof}
By arguments similar to the ones used in the first step of the proof of Theorem~\ref{th:mainBer}, we have
\[
\bigl\{ A_{t+1} = a \bigr\} \subseteq
\Bigl\{ \mu^\star-\epsilon < \hat{\mu}_{a,N_t(a)} \Bigr\} \, \cup \,
\Bigl\{ \mu^\star-\epsilon > B^+_{a^\star,t} \Bigr\} \, \cup \,
\Bigl\{ \mu^\star-\epsilon \in \bigl[ \hat{\mu}_{a,N_t(a)}, \,\, B^+_{a,t} \bigr] \Bigr\}\,;
\]
indeed, on the event
\qquad
$\displaystyle{\bigl\{ A_{t+1} = a \bigr\} \, \cap \,
\Bigl\{ \mu^\star-\epsilon \geq \hat{\mu}_{a,N_t(a)} \Bigr\} \, \cap \,
\Bigl\{ \mu^\star-\epsilon \leq B^+_{a^\star,t} \Bigr\}}\,,
$ \\
we have, $\hat{\mu}_{a,N_t(a)} \leq \mu^\star-\epsilon \leq B^+_{a^\star,t} \leq B^+_{a,t}$ (where the last
inequality is by definition of the strategy).
Before proceeding, we note that
\[
\Bigl\{ \mu^\star-\epsilon \in \bigl[ \hat{\mu}_{a,N_t(a)}, \,\, B^+_{a,t} \bigr] \Bigr\}
\subseteq \Bigl\{ N_t(a) \,\, \kmin \bigl( \hat{\nu}_{a,N_t(a)}, \, \mu^\star-\epsilon \bigr) \leq f(t) \Bigr\}\,,
\]
since $\kmin$ is a non-decreasing function in its second argument and $\kmin\bigl(\nu,E(\nu)\bigr) = 0$
for all distributions $\nu$.
Therefore,
\begin{multline}
\nonumber
\E \bigl[ N_T(a) \bigr] \leq 1 + \sum_{t=|\cA|}^{T-1} \P \Bigl\{ \mu^\star-\epsilon < \hat{\mu}_{a,N_t(a)} \andAt \Bigr\}
+ \sum_{t=|\cA|}^{T-1} \P \Bigl\{ \mu^\star-\epsilon > B^+_{a^\star,t} \Bigr\} \\
+ \sum_{t=|\cA|}^{T-1} \P \Bigl\{ N_t(a) \,\, \kmin \bigl( \hat{\nu}_{a,N_t(a)}, \, \mu^\star-\epsilon \bigr) \leq f(t) \andAt \Bigr\}\,;
\end{multline}
now, the two sums with the events ``and $A_{t+1} = a$'' can be rewritten by using the stopping times $\tau_{a,k}$
introduced in the proof of Theorem~\ref{th:mainBer}; more precisely, by mimicking the transformations
performed in its step 3, we get the simpler bound
\begin{multline}
\label{eq:threesums}
\E \bigl[ N_T(a) \bigr] \leq 1 + \sum_{k=2}^{T-|\cA|+1} \P \Bigl\{ \mu^\star-\epsilon < \wt{\mu}_{a,k-1} \Bigr\}
+ \sum_{t=|\cA|}^{T-1} \P \Bigl\{ \mu^\star-\epsilon > B^+_{a^\star,t} \Bigr\} \\
+ \sum_{k=2}^{T-|\cA|+1} \P \Bigl\{ (k-1) \,\, \kmin \bigl( \wt{\nu}_{a,k-1}, \, \mu^\star-\epsilon \bigr) \leq f(t) \Bigr\}\,,
\end{multline}
where the $\wt{\nu}_{a,s}$ and $\wt{\mu}_{a,s}$ are respectively the empirical distributions and empirical
expectations computed on the first $s$ elements of the sequence of the random variables $\wt{X}_{a,j} = Y_{\tau_{a,j}}$,
which are {i.i.d.} according to $\nu_a$.

\textbf{Step 1: The first sum in~(\ref{eq:threesums})} is bounded by resorting to Hoeffding's inequality,
whose application is legitimate since $\mu^\star-\mu_a-\epsilon > 0$;
\beqan
\nonumber
\sum_{k=2}^{T-|\cA|+1} \P \Bigl\{ \mu^\star-\epsilon < \wt{\mu}_{a,k-1} \Bigr\}
&=& \sum_{k=1}^{T-|\cA|} \P \Bigl\{ \mu^\star-\mu_a-\epsilon < \wt{\mu}_{a,k} - \mu_a \Bigr\} \\
&\leq& \sum_{k=1}^{T-|\cA|} e^{- 2 k (\mu^\star-\mu_a-\epsilon)^2} \leq \frac{1}{1 - e^{- 2 (\mu^\star-\mu_a-\epsilon)^2}}
\leq \frac{1}{(\mu^\star-\mu_a-\epsilon)^2}\,,
\eeqan
where we used for the last inequality the general upper bounds provided at the beginning of step~5 in the proof of
Theorem~\ref{th:mainBer}.

\textbf{Step 2: The second sum in~(\ref{eq:threesums})} is bounded by first
using the definition of $B^+_{a^\star,t}$, then, decomposing the event depending on the
values taken by $N_t(a^\star)$; and finally using the fact that on $\bigl\{ N_t(a^\star) = k \bigr\}$,
we have the rewriting
\ $
\hat{\nu}_{a,N_t(a)} = \widetilde{\nu}_{a,k}
$ \ and
$ \
\hat{\mu}_{a,N_t(a)} = \widetilde{\mu}_{a,k}\,;
$ \
more precisely,
\beqan
\sum_{t=|\cA|}^{T-1} \P \Bigl\{ \mu^\star-\epsilon > B^+_{a^\star,t} \Bigr\}
 &\leq&  \sum_{t=|\cA|}^{T-1} \P \Bigl\{ N_t(a^\star) \,\, \kmin \bigl( \hat{\nu}_{a^\star,N_t(a^\star)}, \,
\mu^\star-\epsilon \bigr) > f(t) \Bigr\} \\
 &=& \sum_{t=|\cA|}^{T-1} \sum_{k=1}^t \P \Bigl\{ N_t(a^\star) = k \ \,\, \mbox{\small and} \ \,\, k \,\,
\kmin \bigl( \wt{\nu}_{a^\star,k}, \, \mu^\star-\epsilon \bigr) > f(t) \Bigr\} \\
 &\leq&  \sum_{k=1}^T \sum_{t=|\cA|}^{T-1} \P \Bigl\{ k\,\,\kmin \bigl( \wt{\nu}_{a^\star,k}, \, \mu^\star-\epsilon \bigr) > f(t) \Bigr\}\,.
\eeqan
Since $f = \log$ is increasing, we can rewrite the bound, using a Fubini-Tonelli argument, as
\beqan
\sum_{t=|\cA|}^{T-1} \P \Bigl\{ \mu^\star-\epsilon > B^+_{a^\star,t} \Bigr\}
& \leq & \sum_{k=1}^T \,\, \sum_{t=|\cA|}^{T-1} \P \biggl\{ f^{-1} \Bigl( k \, \kmin \bigl( \wt{\nu}_{a^\star,k}, \, \mu^\star-\epsilon \bigr) \Bigr) > t \biggr\} \\
& \leq & \sum_{k=1}^T \, \E \biggl[ f^{-1} \Bigl(k \, \kmin\bigl( \wt{\nu}_{a^\star,k}, \, \mu^{\star}-\epsilon \bigr) \Bigr)
\,\, \mathbb{I}_{ \bigl\{ \kmin( \wt{\nu}_{a^\star,k}, \, \mu^{\star}-\epsilon ) > 0 \bigr\} } \biggr]\,.
\eeqan

Now, \citet[Lemma~13]{Honda10} indicates that, since $\mu^\star - \epsilon \in [0,1)$,
\[
\sup_{\nu \in \cP_F ( [0,1] )} \kmin\bigl(\nu,\mu^\star - \epsilon\bigr) \leq \log \bigl( 1/ (1-\mu^\star + \epsilon) \bigr) \eqdef K_{\max}\,;
\]
we define $Q = K_{\max}/\epsilon^2$ and introduce the following sets $(V_q)_{1 \leq q\leq Q}$:
\beqan
V_q = \Bigl\{ \nu \in \cP_F \bigl( [0,1] \bigr) : \ \ (q-1)\epsilon^2 < \kmin\bigl(\nu,\mu^*-\epsilon) \leq q \epsilon^2 \Bigr\}.
\eeqan

A peeling argument (and by using that $f^{-1} = \exp$ is increasing as well) entails, for all $k \geq 1$,
\beqa
\lefteqn{
\E \biggl[ f^{-1} \Bigl(k \, \kmin\bigl( \wt{\nu}_{a^\star,k}, \, \mu^{\star}-\epsilon \bigr) \Bigr)
\,\, \mathbb{I}_{ \bigl\{ \kmin( \wt{\nu}_{a^\star,k}, \, \mu^{\star}-\epsilon ) > 0 \bigr\} } \biggr] } \\
& = & \sum_{q=1}^Q \, \E \biggl[ f^{-1} \Bigl(k \, \kmin\bigl( \wt{\nu}_{a^\star,k}, \, \mu^{\star}-\epsilon \bigr) \Bigr)
\,\, \mathbb{I}_{ \bigl\{ \wt{\nu}_{a^\star,k} \in V_q \bigr\} } \biggr]
\notag\\
& \leq & \sum_{q=1}^Q \, \P \bigl\{ \wt{\nu}_{a^\star,k} \in V_q \bigr\} \, f^{-1}(kq\epsilon^2)
\leq
\sum_{q=1}^Q \P \Bigl\{ \kmin\bigl( \wt{\nu}_{a^\star,k}, \, \mu^{\star} -\epsilon \bigr) > (q-1)\epsilon^2 \Bigr\} \, f^{-1}(kq\epsilon^2)\,,
\label{eq:peeling}
\eeqa
where we used the definition of $V_q$ to obtain each of the two inequalities.
Now, by Lemma~\ref{lm:dminvar}, when $E \bigl( \wt{\nu}_{a^\star,k} \bigr) < \mu^\star - \epsilon$,
which is satisfied whenever $\kmin \bigl( \wt{\nu}_{a^\star,k}, \, \mu^\star-\epsilon \bigr) > 0$,
we have
\[
\kmin \bigl( \wt{\nu}_{a^\star,k}, \, \mu^\star-\epsilon \bigr)
\leq \kmin \bigl( \wt{\nu}_{a^\star,k}, \, \mu^\star \bigr) - 2 \epsilon^2
\leq \cK \bigl( \wt{\nu}_{a^\star,k}, \, \nu^\star \bigr) - 2 \epsilon^2 \,,
\]
where the last inequality is by mere definition of $\kmin$.
Therefore,
\[
\P \Bigl\{ \kmin \bigl( \wt{\nu}_{a^\star,k}, \, \mu^\star-\epsilon \bigr) > (q-1)\epsilon^2 \Bigr\}
\leq \P \Bigl\{ \cK \bigl( \wt{\nu}_{a^\star,k}, \, \nu^\star \bigr) >  (q+1)\epsilon^2 \Bigr\}\,.
\]

We note that for all $k \geq 1$, \qquad
$\displaystyle{
\P \Bigl\{ \cK \bigl( \wt{\nu}_{a^\star,k}, \, \nu^\star \bigr) > (q+1)\epsilon^2 \Bigr\}
\leq (k+1)^{|\cS^\star|} \, e^{-k(q+1)\epsilon^2}\,,
}$ \\
where we recall that $\cS^\star$ denotes the finite support of $\nu^\star$ and
where we applied \ifarXiv{Corollary~\ref{cor:types} of the appendix.}\ifCOLT{the method
of types; see, e.g., the extended version of the present paper~\citep{Ext} for more
details about this standard inequality.}
Now, (\ref{eq:peeling}) then yields, via the choice
$f = \log$ and thus $f^{-1} = \exp$, that
\[
\E \biggl[ f^{-1} \Bigl(k \, \kmin\bigl( \wt{\nu}_{a^\star,k}, \, \mu^{\star}-\epsilon \bigr) \Bigr)
\,\, \mathbb{I}_{ \bigl\{ \kmin( \wt{\nu}_{a^\star,k}, \, \mu^{\star}-\epsilon ) > 0 \bigr\} } \biggr]
\leq \underbrace{\sum_{q=1}^Q (k+1)^{|\cS^\star|} \, e^{-k(q+1)\epsilon^2} e^{k q \epsilon^2}}_{
= Q \, (k+1)^{|\cS^\star|} \, e^{-k\epsilon^2}}\,.
\]
Substituting the value of $Q$, we therefore have proved that
\[
\sum_{t=|\cA|}^{T-1} \P \Bigl\{ \mu^\star-\epsilon > B^+_{a^\star,t} \Bigr\}
\leq \frac{1}{\epsilon^2}\log \biggl( \frac{1}{1-\mu^* + \epsilon} \biggr) \sum_{k=1}^T (k+1)^{|\cS^\star|} \, e^{-k\epsilon^2}.
\]

\textbf{Step 3: The third sum in~(\ref{eq:threesums})} is first upper bounded by Lemma~\ref{lm:dminvar}, which states that
\[
\kmin \bigl( \wt{\nu}_{a,k-1}, \, \mu^\star \bigr) - \epsilon/ (1 - \mu^\star)
\leq \kmin \bigl( \wt{\nu}_{a,k-1}, \, \mu^\star-\epsilon \bigr)\,,
\]
for all $k \geq 1$, and by using $f(t) \leq f(T)$; this gives
\begin{multline}
\nonumber
\sum_{k=1}^{T-|\cA|} \P \Bigl\{ k \,\, \kmin \bigl( \wt{\nu}_{a,k}, \, \mu^\star-\epsilon \bigr) \leq f(t) \Bigr\} \\
\leq \sum_{k=1}^{T-|\cA|} \P \left\{ k \,\, \kmin \bigl( \wt{\nu}_{a,k}, \, \mu^\star \bigr) \leq f(T) +
\frac{k \, \epsilon}{1 - \mu^\star} \right\}
= \sum_{k=1}^{T-|\cA|} \P \Bigl\{ \wt{\nu}_{a,k} \in \ol{\cC}_{\mu^\star,\gamma_k} \Bigr\}\,,
\end{multline}
where $\gamma_k = f(T)/k + \epsilon/(1-\mu^\star)$ and where the set $\ol{\cC}_{\mu^\star,\gamma_k}$ was defined
in Section~\ref{sec:setkmin}.
For all $\gamma > 0$, we then introduce
\[
\theta_a(\gamma) = \inf \Bigl\{ \cK(\nu',\nu_a) : \ \ \nu' \in \cC_{\mu^\star,\gamma} \Bigr\}
= \inf \Bigl\{ \cK(\nu',\nu_a) : \ \ \nu' \in \ol{\cC}_{\mu^\star,\gamma} \Bigr\}\,,
\]
(where the second equality follows from the lower semi-continuity of $\cK$)
and aim at bounding $\P \Bigl\{ \wt{\nu}_{a,k} \in \ol{\cC}_{\mu^\star,\gamma} \Bigr\}$.

As shown in Section~\ref{sec:setkmin},
the set $\cC_{\mu^\star,\gamma}$ is a non-empty open convex set.
If we prove that $\theta_a(\gamma)$ is finite for all $\gamma > 0$, then all the conditions will be required to apply
Lemma~\ref{lm:Sanov2} and get the upper bound
\[
\sum_{k=1}^{T-|\cA|} \P \Bigl\{ \wt{\nu}_{a,k} \in \ol{\cC}_{\mu^\star,\gamma_k} \Bigr\}
\leq \sum_{k=1}^{T-|\cA|} \, e^{-k \, \theta_a(\gamma_k)}\,.
\]

To that end, we use the fact that $\nu_a$ is finitely supported.
Now, either the probability of interest is null and we are done; or, it is not null, which implies
that there exists a possible value of $\wt{\nu}_{a,k}$ that is in $\ol{\cC}_{\mu^\star,\gamma}$; since
this value is a distribution with a support included in the one of $\nu_a$, it is absolutely continuous
with respect to $\nu_a$ and hence, the Kullback-Leibler divergence between this value and $\nu_a$
is finite; in particular, $\theta_a(\gamma)$ is finite.

Finally, we bound the $\theta_a(\gamma_k)$ for values of $k$ larger than \qquad
$\displaystyle{
k_0 = \left\lceil \frac{(1+c_a)\,f(T)}{\kmin(\nu_a,\mu^\star)} \right\rceil\,;
}$ \\
we have that for all $k \geq k_0$, in view of the bound put
on $\epsilon$,
\begin{equation}
\label{eq:gammak}
\gamma_k \leq \gamma_{k_0} = \frac{f(T)}{k_0} + \frac{\epsilon}{1-\mu^\star} < \frac{\kmin(\nu_a,\mu^\star)}{1+c_a} +
\frac{c_a/2}{1+c_a} \, \kmin(\nu_a,\mu^\star) = \frac{1+c_a/2}{1+c_a} \, \kmin(\nu_a,\mu^\star)\,.
\end{equation}
Since $\theta_a$ is non increasing, we have
\[
\sum_{k=1}^{T-|\cA|} \, e^{-k \, \theta_a(\gamma_k)}
\leq k_0-1 + \sum_{k=k_0}^{T-|\cA|} \, e^{-k \, \theta_a(\gamma_{k_0})}
\leq k_0-1 + \frac{1}{1 - e^{- \Theta_a(c_a,\epsilon)}}\,,
\]
provided that the quantity $\Theta_a(c_a,\epsilon) = \theta_a \bigl( \gamma_{k_0} \bigr)$
is positive, which we prove now.

Indeed for all $\nu'\in\cC_{\mu^\star,\gamma_{k_0}}$, we have by definition and by~(\ref{eq:gammak}) that
\[
\kmin(\nu',\mu^\star) - \kmin(\nu_a,\mu^\star) < \gamma_{k_0} - \kmin(\nu_a,\mu^\star) < - \bigl( (c_a/2) \big/ (1+c_a) \bigr)
\kmin(\nu_a,\mu^\star)\,.
\]
Now, in the case where $\E_{\nu_a} \bigl[ (1-\mu^\star)/(1-X) \bigr] > 1$,
we have, first by application of Pinsker's inequality and then by Lemma~\ref{lm:pinsker}, that
\[
\cK\bigl( \nu',\nu_a\bigr) \, \geq  \, \frac{\Arrowvert \nu'-\nu_a \Arrowvert^2_1}{2} \,
\geq \, \frac{1}{2 \, M_{\nu_a,\mu^\star}^2} \bigl( \kmin(\nu_a,\mu^\star) - \kmin(\nu',\mu^\star) \big)^2
> \, \frac{c_a^2 \, \bigl( \kmin(\nu_a,\mu^\star) \bigr)^2}{8 \, (1+c_a)^2 \, M_{\nu_a,\mu^\star}^2} \,;
\]
since, again by Pinsker's inequality, $\kmin(\nu_a,\mu^\star) \geq (\mu_a - \mu^\star)^2/2 > 0$,
we have exhibited a lower bound independent of $\nu'$ in this case.
In the case where $\E_{\nu_a} \bigl[ (1-\mu^\star)/(1-X) \bigr] \leq 1$,
we apply the second part of Lemma~\ref{lm:pinsker}, with $\alpha_a = (c_a/2)/(1+c_a)$, and get
\[
\cK\bigl( \nu',\nu_a\bigr) \, \geq  \, \frac{\Arrowvert \nu'-\nu_a \Arrowvert^2_1}{2} \, \geq\,
\frac{1}{2} \, \left( \frac{1-\mu^\star}{(2/\alpha_a) \, \bigl( (2/\alpha_a) - 1 \bigr)} \right)^2 > 0\,.
\]
Thus, in both cases we found a positive lower bound independent of $\nu'$,
so that the infimum over $\nu' \in \cC_{\mu^\star,\gamma_{k_0}}$ of the quantities
$\kmin(\nu',\mu^\star)$, which precisely equals $\theta_a \bigl( \gamma_{k_0} \bigr)$, is also positive.
This concludes the proof.
\end{proof}

{\small
\paragraph{Conclusion.}
We provided a finite-time analysis of the (asymptotically optimal) $\kmin$--strategy in the case of finitely supported distributions. \ifarXiv{One could think that the extension to the case of general distributions is straightforward. However this extension appears somewhat difficult (at least when using the current definition of $\kmin$) for the following reasons: (1) Step 2 in the proof uses the method of types,
that would require some extension of Sanov's non-asymptotic Theorem to this case. (2) Step 3 requires to have both $\theta_a(\gamma)<\infty$ for all $\gamma>0$ and $\theta_a(\gamma)>0$ for $\gamma < \kmin(\nu_a,\mu^\star)$, which does not seem to be always the case for general distributions.
Exploring other directions for such extensions is left for future work; for instance, histogram-based approximations of general distributions
could be considered.}\ifCOLT{The extension to the case of general distributions (e.g., by histogram-based approximations of
such general distributions) is left for future work.}
}


%
%
%
%

%
%
%
%
%
%

{\small \paragraph{Acknowledgements.}
The authors wish to thank \textit{Peter Auer} and \textit{Daniil Ryabko} for insightful discussions.
They acknowledge support from
the French National Research Agency (ANR)
under grant EXPLO/RA (``Exploration--exploitation for efficient resource allocation'')
and by the PASCAL2 Network of Excellence under EC grant {no.} 506778.
}


\ifarXiv{\newpage}

{\small
\bibliography{Bib-Japonais}
}

\ifarXiv{
\newpage

\appendix
\section{Appendix beyond the COLT page limit}

A conference version of this paper was published in the
\emph{Proceedings of the Twenty-Fourth Annual Conference on Learning Theory} (COLT'11);
this appendix details some material which was alluded at in this conference version
but could not be published therein because of the page limit.

\subsection{Proof of Lemma~\ref{lm:Aure}}

\emph{We only provide it for the convenience of the readers} since it is similar to
the one presented in~\citet[Proposition~2]{GaLe10} or in~\citet{GaCa11};
it was however somewhat simplified by noting that the proof technique used
leads to a maximal inequality, as stated in Lemma~\ref{lm:Aure}, and not only to an inequality
for a self-normalized average, as stated in the original reference. \\

\begin{proof}
The result is straightforward in the cases $p = 0$ or $p = 1$, since then, $\wh{p}_s = p$
almost surely; in the rest of the proof, we therefore only consider the case where $p \in (0,1)$.

It suffices to show the first bound stated in the lemma, since the second one follows
by a decomposition of the probability space according to the values of $N_t$.
Actually, we will show
\[
\P_p \! \left( \bigcup_{s = 1}^t  \biggl\{ s\,\,\cK\Bigl( \beta\bigl({\wh{p}_s}\bigr), \, \beta(p) \Bigr) \geq \eps
\ \,\, \mbox{\small and} \ \,\, \wh{p}_s > p \biggr\} \right) \leq e \, \bigl\lceil \eps \log t \bigr\rceil \, e^{-\eps}\,,
\]
and the desired result will follow by symmetry and a union bound.

\textbf{Step 1: A martingale.} For all $\lambda > 0$, we consider the log-Laplace transform
\[
\psi_p(\lambda) = \log \E_p \bigl[ e^{\lambda X_1} \bigr]
= \log \bigl( (1-p) + p\,e^\lambda \bigr)\,,
\]
with which we define the martingale
\[
W_s(\lambda) = \exp \bigl( \lambda (X_1 + \ldots + X_s) - s \, \psi_p(\lambda) \bigr)\,.
\]

\textbf{Step 2: A peeling argument.} We introduce $t_0 = 1$
and $t_k = \lfloor \gamma^k \rfloor$, for some $\gamma > 1$ that will be defined by
the analysis. We also denote by $K = \bigl\lceil (\log t)/(\log \gamma) \bigr\rceil$
an upper bound on the number of elements in the peeling.

We also note that by continuity of the Kullback-Leibler divergence in the
case of Bernoulli distributions, for all $\epsilon > 0$, there exists a unique
element $p_\epsilon \in (p,1)$ such that $\cK \bigl( \beta({q_\epsilon}), \, \beta(p) \bigr) = \epsilon$;
this element satisfies that
\[
\cK \bigl( \beta(q), \, \beta(p) \bigr) \geq \epsilon \ \ \mbox{and} \ \ q \geq p \qquad \mbox{entails} \qquad
q \geq p_{\epsilon}\,.
\]

Denoting by $\epsilon_k = \epsilon/t_k$, a union bound using the described peeling then yields
\begin{eqnarray*}
\lefteqn{ \P_p \! \left( \bigcup_{s = 1}^t  \biggl\{ s\,\,\cK\Bigl( \beta\bigl({\wh{p}_s}\bigr), \, \beta(p) \Bigr) \geq \eps
\ \,\, \mbox{\small and} \ \,\, \wh{p}_s > p \biggr\} \right) } \\
& \leq & \sum_{k=1}^K \, \P_p \! \left( \bigcup_{s = t_{k-1}}^{t_k}  \biggl\{ s\,\,\cK\Bigl( \beta\bigl({\wh{p}_s}\bigr), \, \beta(p) \Bigr) \geq \eps
\ \,\, \mbox{\small and} \ \,\, \wh{p}_s > p \biggr\} \right) \\
& \leq & \sum_{k=1}^K \, \P_p \! \left( \bigcup_{s = t_{k-1}}^{t_k}  \biggl\{ \cK\Bigl( \beta\bigl({\wh{p}_s}\bigr), \, \beta(p) \Bigr) \geq \frac{\eps}{t_k}
\ \,\, \mbox{\small and} \ \,\, \wh{p}_s > p \biggr\} \right) \\
& = & \sum_{k=1}^K \, \P_p \! \left( \bigcup_{s = t_{k-1}}^{t_k}  \Bigl\{ \wh{p}_s \geq p_{\eps_k} \Bigr\} \right)
\ = \ \sum_{k=1}^K \, \P_p \! \left( \bigcup_{s = t_{k-1}}^{t_k}  \Bigl\{ X_1 + \ldots + X_s - s \, p_{\eps_k} \geq 0 \Bigr\} \right)
\end{eqnarray*}
Now, the variational formula for Kullback-Leibler divergences shows that for all $k$, there exists a $\lambda_k$ such that
\[
\epsilon_k = \cK \bigl( \beta(p_{\epsilon_k}), \, \beta(p) \bigr) = \lambda_k \, p_{\epsilon_k} - \psi_{p}(\lambda_k)\,;
\]
actually, a straightforward calculation shows that $\lambda_k = \log \bigl( p_{\epsilon_k} (1-p \bigr) -
\log \bigl( p (1-p_{\epsilon_k}) \bigr) > 0$ is a suitable value.
Thus,
\begin{eqnarray*}
\lefteqn{ \sum_{k=1}^K \, \P_p \! \left( \bigcup_{s = t_{k-1}}^{t_k}  \Bigl\{ X_1 + \ldots + X_s - s \, p_{\eps_k} \geq 0 \Bigr\} \right) } \\
& = & \sum_{k=1}^K \, \P_p \! \left( \bigcup_{s = t_{k-1}}^{t_k}  \Bigl\{ \exp \bigl( \lambda_k (X_1 + \ldots + X_s) - \lambda_k s \, p_{\eps_k} \bigr) \geq 1 \Bigr\} \right) \\
& = & \sum_{k=1}^K \, \P_p \! \left( \bigcup_{s = t_{k-1}}^{t_k}  \Bigl\{ \exp \bigl( \lambda_k (X_1 + \ldots + X_s) - s \, \psi_{p}(\lambda_k) \bigr) \geq
e^{s \, \epsilon_k}  \Bigr\} \right) \\
& \leq & \sum_{k=1}^K \, \P_p \! \left( \bigcup_{s = t_{k-1}}^{t_k}  \Bigl\{ W_s(\lambda_k) \geq e^{t_{k-1} \epsilon_k}  \Bigr\} \right) \\
& \leq & \sum_{k=1}^K \, e^{-t_{k-1} \, \epsilon_k} = K e^{-\epsilon/\gamma}\,,
\end{eqnarray*}
where in the last step, we resorted to Doob's maximal inequality.

\textbf{Step 3: Choosing $\gamma$.} The obtained bound equals, by substituting the value of $K$
and by choosing $\gamma = \epsilon/(\epsilon - 1)$,
\[
K e^{-\epsilon/\gamma} = \bigl\lceil (\log t)/(\log \gamma) \bigr\rceil \, e^{-\epsilon + 1}
= \left\lceil \frac{\log t}{\log \bigl( \epsilon /(\epsilon - 1) \bigr)} \right\rceil \, e^{-\epsilon + 1}\,;
\]
the proof is concluded by noting that $\epsilon > 1 \longmapsto \log \bigl( \epsilon /(\epsilon - 1) \bigr) - 1/\epsilon$
is decreasing (its derivative is negative), with limit $0$ at $+\infty$.
\end{proof}

\subsection{Details of the adaptation leading to Lemma~\ref{lm:Sanov2}}

The exact statement of \citet[Theorem~2.1 and comments on page~372]{Din92}
is the following.

\begin{lemma}[Non-asymptotic Sanov's lemma]
\label{lm:Sanov}
Let $\cC$ be an open convex subset of $\cP(\cX)$ such that
\[
\Lambda(\cC) = \inf_{\nup \in \cC} \, \cK(\nup,\nu) < \infty\,.
\]
Then, for all $t \geq 1$,
\[
\P_\nu \bigl\{ \hat{\nu}_t \in \cC \bigr\} \leq e^{- t \Lambda(\ol{\cC})}\,.
\]
\end{lemma}

We show how it entails Lemma~\ref{lm:Sanov2}.
Let $\cC$ be an open convex subset of $\cP(\cX)$
and let $\ol{\cC}$ be its closure.
We denote by
\[
\cC_\delta = \bigl\{ \nu \in \cC : \ \ d(\nu,\cC) < \delta \bigr\}
\]
the $\delta$--open neighborhood of $\cC$, we have $\ol{\cC} \subseteq \cC_\delta$ for all $\delta > 0$.
Therefore, by the lemma above, since $\Lambda(\cC_\delta) \leq \Lambda(\cC) < \infty$,
\[
\P_\nu \bigl\{ \hat{\nu}_t \in \ol{\cC} \bigr\}
\leq \P_\nu \bigl\{ \hat{\nu}_t \in \cC_\delta \bigr\} \leq e^{- t \Lambda(\cC_\delta)}\,.
\]
We pick for each integer $n \geq 1$ an element $\nup_n$ such that $\Lambda\bigl(\cC_{1/n}\bigr) = \cK(\nup_n,\nu) - 1/n$;
by \citet[proof of Proposition 1.1]{Din92}, the sequence of the $\nup_n$ admits a converging subsequence $\nup_{\varphi(n)}$,
whose limit point $\nup_\infty$ belongs to $\ol{\cC}$ and which satisfies
\[
\cK(\nup_\infty,\nu) \leq \liminf_{n \to \infty} \cK(\nup_n,\nu) =
\liminf_{\delta \to 0} \Lambda\bigl(\cC_{\delta}\bigr)\,.
\]
Therefore, by taking limits in the above inequality, we have proved the desired inequality,
\[
\P_\nu \bigl\{ \hat{\nu}_t \in \ol{\cC} \bigr\} \leq e^{- t \cK(\nup_\infty,\nu)} \leq e^{- t \Lambda(\ol{\cC})}\,.
\]

\subsection{Useful properties of $\kmin$ and its level sets}

\paragraph{Proof of Lemma~\ref{lm:pinsker}:}
We resort to the formulation of $\kmin$ in terms of a convex optimization problem as introduced in \cite{HondaArxiv}; more precisely,
it is shown therein that
\begin{equation}
\label{eq:Hvar}
\kmin(\nu,\mu) = \max \biggl\{  \E_\nu \Bigl[ \log \bigl( 1+\lambda(\mu-X) \bigr) \Bigr] : \ \  \lambda \in \bigl[0, \, 1/(1-\mu) \bigr] \biggr\}
\end{equation}
(where $X$ denotes a random variable distributed according to $\nu$), as well as the following
alternative.
The optimal value $\lambda_\nu$ of the parameter $\lambda$ indexing the set is
equal to $1/(1-\mu)$ if and only if $\E_\nu \bigl[ (1-\mu)/(1-X) \bigr] \leq 1$,
and lies in $\bigl[ 0, \, 1/(1-\mu) \bigr)$ if $\E_\nu \bigl[ (1-\mu)/(1-X) \bigr] > 1$.

For all $\lambda \in \bigl[0, \, 1/(1-\mu) \bigr]$, we now introduce the function
\[
\phi_\lambda : x \in [0,1] \,\, \longmapsto \,\, \log \bigl( 1+\lambda(\mu-x) \bigr)\,,
\]
which is always continuous on $[0,1)$; we note also that it is continuous and finite at $x = 1$ when $\lambda < 1/(1-\mu)$.
In the latter case, $\phi_\lambda$ is bounded; since it is decreasing, it is easy to get a uniform bound: for all $x$,
\[
\bigl| \phi_\lambda(x) \bigr| \leq \bigl|\phi(0)\bigr| + \bigl| \phi(1) \bigr| = \log \frac{1+\lambda\mu}{1+\lambda(\mu-1)} \eqdef M_\lambda\,.
\]
It then follows that for all $\lambda \in \bigl[0, \, 1/(1-\mu) \bigr)$,
\begin{equation}
\label{eq:ctrlam}
\E_\nu \bigl[ \phi_{\lambda}(X) \bigr] - \E_{\nu'} \bigl[ \phi_{\lambda}(X) \bigr]
\leq M_{\lambda} \, \Arrowvert \nu - \nu' \Arrowvert_1\,.
\end{equation}

In the case when $\lambda_\nu < 1/(1-\mu)$, we have from the variational formulation~(\ref{eq:Hvar}) that
\[
\kmin(\nu,\mu) - \kmin(\nu',\mu) \leq \E_\nu \bigl[ \phi_{\lambda_\nu}(X) \bigr] - \E_{\nu'} \bigl[ \phi_{\lambda_\nu}(X) \bigr]
\leq M_{\lambda_\nu} \, \Arrowvert \nu - \nu' \Arrowvert_1\,.
\]
Thus, the constant $M_{\nu,\mu}$ in the statement of the lemma
corresponds to our quantity $M_{\lambda_\nu}$ in this case.

We now consider the case where $\lambda_\nu = 1/(1-\mu)$.
By~(\ref{eq:ctrlam}) and variational formulation~(\ref{eq:Hvar}), we have that for all $\lambda \in \bigl[0, \, 1/(1-\mu) \bigr)$,
\begin{multline}
\nonumber
\kmin(\nu,\mu) - \kmin(\nu',\mu)
\leq \kmin(\nu,\mu) - \E_{\nu'} \bigl[ \phi_{\lambda}(X) \bigr] \\
= \Bigl( \kmin(\nu,\mu) - \E_\nu \bigl[ \phi_{\lambda}(X) \bigr] \Bigr) +
\Bigl( \E_\nu \bigl[ \phi_{\lambda}(X) \bigr] - \E_{\nu'} \bigl[ \phi_{\lambda}(X) \bigr] \Bigr)\,.
\end{multline}
The second difference is bounded according to~(\ref{eq:ctrlam}); the first difference is bounded by
concavity of $\lambda < 1/(1-\mu) \,\, \mapsto \,\, \phi_\lambda(x)$, for all $x$:
\begin{multline}
\nonumber
\E_\nu \bigl[ \phi_{\lambda}(X) \bigr] \geq
\bigl(1 - \lambda(1-\mu) \bigr) \, \E_\nu \bigl[ \phi_{0}(X) \bigr]
+ \lambda(1-\mu) \, \E_\nu \bigl[ \phi_{0}(X) \bigr] \\
= \lambda(1-\mu) \, \E_\nu \bigl[ \phi_{1/(1-\mu)}(X) \bigr]
= \lambda(1-\mu) \, \kmin(\nu,\mu)\,,
\end{multline}
since $\phi_0$ is the null function and $\lambda_\nu = 1/(1-\mu)$. Putting all pieces together, we have proved that
for all $\lambda \in \bigl[0, \, 1/(1-\mu) \bigr)$,
\begin{equation}
\label{eq:pfin}
\kmin(\nu,\mu) - \kmin(\nu',\mu) \leq \bigl( 1 - \lambda(1-\mu) \bigr) \, \kmin(\nu,\mu)
+ M_{\lambda} \, \Arrowvert \nu - \nu' \Arrowvert_1\,.
\end{equation}
We recall that by assumption,
$\kmin(\nu,\mu) - \kmin(\nu',\mu) \geq \alpha \, \kmin(\nu,\mu)$ with $\alpha \in (0,1)$,
so that the choice $\lambda = (1-\alpha/2)/(1-\mu)$, which indeed lies in $\bigl( 0, \, 1/(1-\mu) \bigr)$,
is such that
\[
M_\lambda = \log \! \left( 1 + \frac{\lambda}{1+\lambda(\mu-1)} \right) = \log \!\left( 1 + \frac{\lambda}{\alpha/2} \right)
\leq \frac{2\lambda}{\alpha}\,,
\]
so that~(\ref{eq:pfin}) entails
\[
\alpha \, \kmin(\nu,\mu) \leq \frac{\alpha}{2} \, \kmin(\nu,\mu) + \frac{2\lambda}{\alpha} \, \Arrowvert \nu - \nu' \Arrowvert_1\,,
\]
and finally
\[
\Arrowvert \nu - \nu' \Arrowvert_1 \geq \frac{\alpha^2}{4 \lambda} = \frac{\alpha^2 (1-\mu)}{1 - \alpha/2} = \frac{1-\mu}{(2/\alpha) \, \bigl( (2/\alpha) - 1 \bigr)}\,;
\]
which concludes the proof.
\qed

\paragraph{Proof of Lemma~\ref{lm:dminvar}:}
In \cite{HondaArxiv} it is shown that in this case, $\kmin(\nu,\mu)$ is differentiable in $\mu \in (E(\nu),1)$ with
\beqa
\frac{1}{1-\mu}\geq \frac{\partial}{\partial \mu}\kmin(\nu,\mu) \geq \frac{\mu- E(\nu)}{\mu(1-\mu)}. \label{eq:partial}
\eeqa
We apply this result to the rewriting
\beqan
\kmin( \nu,\mu) - \kmin( \nu,\mu - \epsilon) = \int_{\mu-\epsilon}^\mu \frac{\partial}{\partial \mu}\kmin(\nu,u) \,\d u \,,
\eeqan
which already gives one part of the bound. For the lower bound, we note that
by assumption $- E(\nu) > - (\mu-\epsilon)$ and that $u(1-u) \leq 1/4$ (since we consider distributions with support included in $[0,1]$);
so that, for all $u$,
\beqan
\frac{u- E(\nu)}{u(1-u)} \geq 4 \bigl( u-(\mu-\epsilon) \bigr) \, .
\eeqan
Integrating the bound concludes the main part of the proof.

Now, to see that the first inequality in the statement is always valid, we need to consider the
case when $E(\nu)\geq\mu$, for which the statement is trivial since then $\kmin(\nu,\mu) = 0$, and the case when $\mu > E(\nu) \geq \mu-\epsilon$.
But in the latter case, it is shown in \citet[Lemma~6, case 2]{HondaArxiv} that
\beqan
\kmin(\nu,\mu) \leq \frac{\mu - E(\nu)}{1-\mu}\,,
\eeqan
which concludes the proof.
\qed

\paragraph{Proof of Lemma~\ref{lm:closure}:}
First, $\cC_{\mu}(\gamma)$ is non empty as it always contains $\delta_{\mu}$, the Dirac mass on $\mu$.

The fact that $\cC_{\mu}(\gamma)$ is convex follows from the convexity of $\cK$ in the pair
of probability distributions that it takes as an argument. Indeed, for all $\alpha \in [0,1]$, $\nu', \, \nu'' \in \cC_{\mu}(\gamma)$, denoting by $\nu'_\mu,\,\nu''_\mu$ some distributions such that
the defining conditions in $\cC_{\mu}(\gamma)$ are satisfied, we have that
\[
E \bigl( \alpha \nu'_\mu + (1-\alpha)\nu''_\mu \bigr) > \mu
\]
and
\[
\cK\bigl( \alpha \nu' + (1-\alpha)\nu'', \, \alpha \nu'_\mu + (1-\alpha)\nu''_\mu \bigr)
\leq \alpha \, \cK\bigl( \nu',\nu'_\mu\bigr) + (1-\alpha) \, \cK\bigl( \nu'',\nu''_\mu\bigr) < \gamma\,.
\]

We prove that $\cC_{\mu}(\gamma)$ is an open set.
With each $\nu' \in \cC_{\mu}(\gamma)$, we associate a distribution
$\nu'_\mu$ satisfying the defining constraints in $\cC_{\mu}(\gamma)$; by choosing
\[
\alpha = \frac{1-\mu \big/ E\bigl(\nu'_\mu\bigr)}{2} \,\, \in (0,\,1/2),
\]
we have that the open set formed by the
\[
(1-\alpha) \, \nu' + \alpha \, \nu'', \qquad \nu'' \in \mbox{B}(\nu',1)
\]
 is contained in $\cC_{\mu,\gamma}$, where $\mbox{B}(\nu',1)$ denotes the ball with center $\nu'$ and radius 1
in the norm $\norm$ over $\cP(\cX)$. Indeed, we have on the one hand,
\[
E\bigl( (1-\alpha) \, \nu'_\mu + \alpha \, \nu'' \bigr) \geq (1-\alpha) \, E\bigl(\nu'_\mu\bigr) \geq \left( 1 - \frac{1-\mu \big/
E\bigl(\nu'_\mu \bigr)}{2} \right) E\bigl(\nu'_\mu\bigr) = \frac{E\bigl(\nu'_\mu\bigr) + \mu}{2} > \mu\,,
\]
and on the other hand, by convexity of the Kullback-Leibler divergence,
\[
\cK \bigl( (1-\alpha) \, \nu' + \alpha \, \nu'', \, (1-\alpha) \, \nu'_\mu + \alpha \, \nu'' \bigr)
\leq (1-\alpha) \, \cK \bigl( \nu', \, \nu'_\mu \bigr) <  (1-\alpha) \gamma\,.
\]

To prove the desired inclusion, we first note that in the case of $\cP_F \bigl( [0,1] \bigr)$,
\cite{HondaArxiv} show that one has the rewriting
\[
\kmin(\nu,\mu) = \min \, \Bigl\{ \cK(\nu,\nu') : \ \ \nu' \in \cP_F \bigl( [0,1] \bigr) \ \ \mbox{\rm s.t.} \ \ E(\nu') \geq \mu \Bigr\}\,;
\]
in particular, the infimum is achieved with this new formulation.
Hence,
\[
\cC_{\mu,\gamma} = \Bigl\{ \nu' \in \cP_F \bigl( [0,1] \bigr) : \ \ \exists \, \nu'_\mu \in \cP_F \bigl( [0,1] \bigr)
\ \ \mbox{s.t.} \ \ E\bigl(\nu'_\mu\bigr) \geq \mu \ \ \mbox{and} \ \ \cK\bigl( \nu',\nu'_\mu\bigr) < \gamma \Bigr\}\,.
\]
Also, an element of the set of interest is therefore a $\nu' \in \cP_F \bigl( [0,1] \bigr)$ such that $\kmin(\nu',\mu) \leq \gamma$,
that is, such that there exists $\nu'_\mu \in \cP \bigl( [0,1] \bigr)$
with $E\bigl(\nu'_\mu\bigr) \geq \mu$ and $\cK\bigl( \nu',\nu'_\mu\bigr) \leq \gamma$.
Now, the distributions
\[
\nu'_n = \left( 1 - \frac{1}{n} \right) \nu' + \frac{1}{n} \delta_1\,, \qquad \mbox{thanks to the} \qquad \nu'_{\mu,n} =
\left( 1 - \frac{1}{n} \right) \nu'_\mu + \frac{1}{n} \delta_1\,,
\]
all belong to $\cC_\gamma$, as, similarly to the above argument,
\[
E \bigl( \nu'_n \bigr) \geq \mu + \frac{1-\mu}{n} > \mu \qquad \mbox{and} \qquad
\cK \bigl( \nu'_n, \, \nu'_{\mu,n} \bigr) \leq \left( 1 - \frac{1}{n} \right) \cK\bigl( \nu',\nu'_\mu\bigr) < \gamma\,.
\]
In addition, we have by construction that the $\nu'_n$ converge to $\nu'$, hence, $\nu' \in \ol{\cC}_\gamma$.
\qed

\subsection{The method of types}

Let $X_1,X_2,\ldots$ be a sequence of random variables that are {i.i.d.} according to
a distribution denoted by $\nu$. In this subsection, we will index all probabilities and
expectations by $\nu$.

For all $k \geq $, we denote by $\cE_k$ the set of
possible values (the so-called types) of the empirical distribution
\[
\hat{\nu}_k = \sum_{j = 1}^k \delta_{X_j}\,.
\]
If $\nu$ has a finite support denoted by $\cS$, then the cardinality $|\cE_k|$ of $\cE_k$ is bounded by $(k+1)^{|\cS|}$.

\begin{lemma}
In the case where $\nu$ has a finite support,
for all $k \geq 1$ and $\nup \in \cE_k$,
\[
\P_\nu \bigl\{ \hat{\nu}_k = \nup \bigr\} \leq e^{-k \, \cK(\nup,\nu)}\,.
\]
\end{lemma}

\begin{corollary}
\label{cor:types}
In the case where $\nu$ has a finite support,
for all $k \geq 1$, all $\gamma > 0$,
\begin{multline}
\nonumber
\P \Bigl\{ \cK \bigl( \hat{\nu}_k, \, \nu \bigr) > \gamma \Bigr\}
= \sum_{\nup \in \cE_k} \ind_{ \{ \cK(\nup,\nu) > \gamma \} } \, \P_\nu \bigl\{ \hat{\nu}_k = \nup \bigr\} \\
\leq \sum_{\nup \in \cE_k} \ind_{ \{ \cK(\nup,\nu) > \gamma \} } \, e^{-k \, \cK(\nup,\nu)}
\leq |\cE_k| \, e^{-k\gamma} \leq (k+1)^{|\cS|} \, e^{-k\gamma} \,.
\end{multline}
\end{corollary}

}

\end{document}